\pdfoutput=1
\RequirePackage{ifpdf}
\ifpdf
\documentclass[pdftex]{sigma}
\else
\documentclass{sigma}
\fi

\usepackage{bbm}

\begin{document}

\allowdisplaybreaks

\renewcommand{\thefootnote}{$\star$}

\renewcommand{\PaperNumber}{038}

\FirstPageHeading

\ShortArticleName{Relative Critical Points}

\ArticleName{Relative Critical Points\footnote{This paper is a~contribution to the Special Issue  ``Symmetries of Dif\/ferential
Equations:  Frames, Invariants and~Applications''.
The full collection is available at
\href{http://www.emis.de/journals/SIGMA/SDE2012.html}{http://www.emis.de/journals/SIGMA/SDE2012.html}}}

\Author{Debra LEWIS}

\AuthorNameForHeading{D.~Lewis}

\Address{Mathematics Department, University of California, Santa Cruz, Santa Cruz, CA 95064, USA}
\Email{\href{mailto:lewis@ucsc.edu}{lewis@ucsc.edu}}
\URLaddress{\url{http://people.ucsc.edu/~lewis/}}

\ArticleDates{Received October 01, 2012, in f\/inal form May 06, 2013; Published online May 17, 2013}

\Abstract{Relative equilibria of Lagrangian and Hamiltonian systems with symmetry are critical points of appropriate scalar functions parametrized by the Lie algebra (or its dual) of the symmetry group. Setting aside the structures~-- symplectic, Poisson, or variational~-- generating dynamical systems from such functions highlights the common features of their construction and analysis, and supports the construction of analogous functions in non-Hamiltonian settings. If the symmetry group is nonabelian, the
functions are invariant only with respect to the isotropy subgroup of the given parameter value. Replacing the parametrized family of functions with a single function on the product
manifold and extending the action using the (co)adjoint action on the algebra or its dual yields a fully invariant  function.
An invariant map can be used to reverse the usual perspective: rather than selecting a parametrized family of functions and f\/inding
their critical points, conditions under which functions will be critical on specif\/ic orbits, typically distinguished by isotropy class,
can be derived. This strategy is illustrated using several well-known mecha\-ni\-cal systems~-- the Lagrange top, the double spherical
pendulum, the free rigid  body, and the Riemann ellipsoids~-- and generalizations of these systems.
}


\Classification{37J15; 53D20; 58E09; 70H33}

\newcommand{\ddt}[1]{{\smallfrac{d#1}{d t}}}
\newcommand{\bigfrac}[2]{{\frac {#1}{#2}}}
\newcommand{\smallfrac}[2]{{\textstyle \frac {#1}{#2}}}
\newcommand{\norm}[1]{\left \Vert #1 \right \Vert}
\newcommand{\diag}[1]{\lsb #1 \rsb}
\newcommand{\lp}{\left ( }
\newcommand{\rp}{\right ) }
\newcommand{\la}{\left \langle}
\newcommand{\ra}{\right \rangle}
\newcommand{\lcb}{\left \{ }
\newcommand{\rcb}{\right \} }
\newcommand{\lsb}{\left [}
\newcommand{\rsb}{\right ]}
\newcommand{\fb}{{\mathfrak b}}
\newcommand{\fg}{{\mathfrak g}}
\newcommand{\fm}{{\mathfrak m}}
\newcommand{\fn}{{\mathfrak n}}
\newcommand{\fq}{{\mathfrak q}}
\newcommand{\fs}{{\mathfrak s}}
\newcommand{\bi}{{\mathbf I}}
\newcommand{\bp}{{\mathbf p}}
\newcommand{\bv}{{\mathbf v}}
\newcommand{\I}{{\mathbb I}}
\newcommand{\id}{{\mathbb 1}}
\newcommand{\bE}{{\mathbb E}}
\newcommand{\FL}{{\mathbb F}L}
\newcommand{\bI}{{\mathbb I}}
\newcommand{\bN}{{\mathbb N}}
\newcommand{\bZ}{{\mathbb Z}}
\newcommand{\dt}{\Delta t}
\newcommand{\be}{{\boldsymbol{e}}}
\newcommand{\br}{{\boldsymbol{r}}}
\newcommand{\bs}{{\boldsymbol{s}}}
\newcommand{\bt}{{\boldsymbol{t}}}
\newcommand{\abe}{\acute {\boldsymbol{e}}}
\newcommand{\tbe}{\grave {\boldsymbol{e}}}
\newcommand{\bg}{{\mathbf g}}
\newcommand{\bm}{{\mathbf m}}
\newcommand{\bq}{{\boldsymbol{q}}}
\newcommand{\bu}{{\mathbf u}}
\newcommand{\boldv}{{\mathbf v}}
\newcommand{\by}{{\mathbf y}}
\newcommand{\bz}{{\mathbf z}}
\newcommand{\J}{{\mathbf J}}
\newcommand{\q}{{\boldsymbol{q}}}
\newcommand{\R}{{\mathbb R}}
\newcommand{\x}{{\mathbf x}}
\newcommand{\ba}{{\boldsymbol{a}}}
\newcommand{\bc}{{\mathbf c}}
\newcommand{\bd}{{\mathbf d}}
\newcommand{\bx}{{\mathbf x}}
\newcommand{\bw}{{\mathbf w}}
\newcommand{\bJ}{{\mathbf J}}
\newcommand{\bC}{{\mathbb C}}
\newcommand{\bbE}{{\mathbb E}}
\newcommand{\bbF}{{\mathbb F}}
\newcommand{\bbI}{{\mathbb I}}
\newcommand{\bP}{{\mathbb P}}
\newcommand{\bR}{{\mathbb R}}
\newcommand{\bbZ}{{\mathbb Z}}
\newcommand{\idm}{{\mathbbm {1}}}
\newcommand{\zero}{{\mathbf 0}}
\newcommand{\cA}{{\mathcal A}}
\newcommand{\cB}{{\mathcal B}}
\newcommand{\cE}{{\mathcal E}}
\newcommand{\cF}{{\mathcal F}}
\newcommand{\cI}{{\mathcal I}}
\newcommand{\cJ}{{\mathcal J}}
\newcommand{\cK}{{\mathcal K}}
\newcommand{\cL}{{\mathcal L}}
\newcommand{\cO}{{\mathcal O}}
\newcommand{\cQ}{{\mathcal Q}}
\newcommand{\cR}{{\mathcal R}}
\newcommand{\cS}{{\mathcal S}}
\newcommand{\cU}{{\mathcal U}}
\newcommand{\cV}{{\mathcal V}}
\newcommand{\cX}{{\mathcal X}}
\newcommand{\half}{{\textstyle \frac 1 2}}
\newcommand{\ad}[1]{{\rm ad}_{#1}}
\newcommand{\ads}[1]{{\rm ad}^*_{#1}}
\newcommand{\Ad}[1]{{\rm Ad}_{#1}}
\newcommand{\Ads}[1]{{\rm Ad}^*_{#1}}
\newcommand{\Adsg}{\mbox{Ad}^*_g}
\newcommand{\adg}{\mbox{ad}_g}
\newcommand{\adsg}{\mbox{ad}^*_g}
\newcommand{\adsh}{\mbox{ad}^*_h}
\newcommand{\rang}[1]{{\rm range}\,[{#1}]}
\newcommand{\rank}[1]{{\rm rank}[{#1}]}
\newcommand{\cay}{\mbox{cay}}
\newcommand{\algo}{F}
\newcommand{\dia}{\mbox{diag}[\ba]}
\newcommand{\eps}{\epsilon}
\newcommand{\dep}[1]{\smallfrac {d \ }{d \epsilon}\left . #1 \right|_{\epsilon = 0}}
\newcommand{\tilf}{\tilde f}
\newcommand{\tG}{\widetilde G}
\newcommand{\Exp}{\mbox{exp}}
\newcommand{\aez}{(A, \bseta, \bszeta)}
\newcommand{\Eta}{N}
\newcommand{\Zeta}{Z}
\newcommand{\bseta}{{\boldsymbol{\eta}}}
\newcommand{\bszeta}{{\boldsymbol{\zeta}}}
\newcommand{\bsxi}{{\boldsymbol{\xi}}}
\newcommand{\bsmu}{{\boldsymbol{\mu}}}
\newcommand{\bsOmega}{{\boldsymbol{\Omega}}}
\newcommand{\bsrho}{{\boldsymbol{\rho}}}
\newcommand{\bssigma}{{\boldsymbol{\sigma}}}
\newcommand{\bzero}{{\boldsymbol{0}}}
\newcommand{\iinfxi}{bI_\bsxi}
\newcommand{\spn}[1]{\mbox{span}\lsb #1 \rsb}
\newcommand{\pd}[2]{\frac {\partial #1}{\partial #2}}
\newcommand{\vint}{\cV_{\rm int}}

\renewcommand{\thefootnote}{\arabic{footnote}}
\setcounter{footnote}{0}

\section{Introduction}

The analysis of steady motions of  dynamical systems with symmetries has played a pivotal role in the
development of geometric mechanics, and the identif\/ication of such motions is central
to the design and interpretation of mathematical models arising in a broad range of disciplines.
Relative equilibria of Hamiltonian or Lagrangian systems with symmetry can be characterized as critical
points of functions, often parametrized by elements of the algebra (or its dual) of the symmetry group.
Steady motions with nontrivial isotropy frequently serve as the starting point for both conceptual and detailed
analyses of systems with symmetry. A variety of methods for the analysis the stability of equilibria and relative
equilibria of conservative systems rely on critical point characterizations of (relative) equilibria
(see, e.g.,~\cite{AbrahamMarsden, GuilStern, HHMM, LMS, LS89b, LS89a,
MarsdenRatiu,MarsdenRatiu86, LMSP89, MarsdenWeinstein, Meyer, Montaldi, PRW04, PRW08, ZombroHolmes}).
Related techniques have been used to study bifurcations of families of (relative) equilibria~-- highly symmetric
states are often more intuitively or analytically accessible than more general motions; once some families of
steady motions are known, a variety of techniques can be applied to identify symmetry-breaking bifurcations from these
families to branches of steady motions that might be dif\/f\/icult to locate directly
(see, e.g., \mbox{\cite{GS87,GSS,  Krupa, Lewis93, Olver86}}).  Many classical and modern
applications, e.g., Lagrange tops and isotropic hyperelasticity, involve both spatial and body symmetries; for such
systems, physically symmetric steady motions such as sleeping tops and steadily rotating spheroidal conf\/igurations,
are readily predicted and visualized.
The relative equilibria for such systems frequently include families of quasi-periodic motions bifurcating
from families with continuous isotropy subgroups. The evolution of the model now known as the Riemann ellipsoids
is a striking example of this process (see Section~\ref{riemann}, \cite{Chandrasekhar, FassoLewis,Riemann, ROSD}).

At present, there appears to be no overall best method for the analysis of relative equilibria~-- each approach performs well
for important classes of problems, but is inapplicable or awkward to implement in many equally important situations. Riemann's
analysis of the existence and stability of steady af\/f\/ine motions of a very simple liquid drop model employs techniques that are
recognizable to the modern reader as precursors of several of the methods described below~\cite{FassoLewis, Riemann}
 (see Section~\ref{riemann});  dif\/ferent strategies
are used for each of the classes of steady motions, exploiting the special features of each class. Thus it seems advantageous to
emphasize the common features of this family of techniques, avoiding wherever possible the invocation of  specialized
structures that are not essential to the steps under consideration. We will focus strictly on the critical points of
the functions related to group actions on manifolds, ignoring any underlying variational, symplectic, or Poisson structure.
Our intent is both to organize a versatile tool kit and to distinguish components of traditional relative equilibrium
analyses that make explicit use of these geometric structures from those that utilize only the group action.

Relative equilibrium methods that make use of an equivariant momentum map typically either explicitly restrict
attention to a single momentum level set in some stages of the analysis or use a Lagrange multiplier approach
to implement a constrained criticality test. If the symmetry group is nonabelian, both of these approaches break
the symmetry of the original system unless the momentum value under consideration has full isotropy, i.e.\ the
isotropy subgroup with respect to the coadjoint action is the entire group. The functions resulting from such approaches
are invariant only with respect to the isotropy subgroup~-- in some situations, e.g.\ symplectic reduction, the full
group does not act on the new manifold of interest; in others, e.g.\ the energy-momentum method, the original
manifold is retained, but the original function is replaced by one that is no longer invariant
with respect to the full action. In the latter situation, `rigid' variations tangent to the group orbit play a central role in
the analysis, but the information captured by such variations cannot be accessed by applying invariant function
theory to the new reduced-symmetry function.

To maintain the full group invariance, we regard the (dual) algebra element as a variable, not a parameter, and extend the
group action to the product manifold using the (co)adjoint action  to obtain a fully invariant function. This allows the use of
invariant function theory for the full symmetry group and facilitates systematic organization of point-algebra element pairs
by isotropy class. We  only require criticality with respect to the original variables. However, utilizing the invariance
of the function, we can replace variations tangent to the group orbit in the manifold with variations tangent to the (co)adjoint
orbit in the algebra, obtaining the so-called rigid criticality conditions.
This strategy can reveal purely symmetry-dependent features of the system.

Steady af\/f\/ine motions of a self-gravitating incompressible f\/luid mass were studied by some of the greatest mathematicians of the
18th and 19th centuries, including Newton, Maclaurin, Jacobi, Dirichlet \cite{Dirichlet}, Dedekind \cite{Dedekind},
and Riemann~\cite{Riemann}; see \cite{Chandrasekhar} for a detailed, albeit politicized, synopsis of the classic literature. The
techniques introduced in the analysis of this system foreshadow many key components of geometric mechanics, bifurcation
theory, and Lie group theory, and are largely driven by symmetry considerations~-- the af\/f\/ine model is equivariant with
respect to actions on ${\rm SL}(3) \times \R^3 \times \R^3$ induced by left and right multiplication on~${\rm SL}(3)$ by~${\rm SO}(3)$ and the
action of $\bbZ_2$ on ${\rm SL}(3)$ by transposition,
as well as an additional~$\bbZ_2^2$ action def\/ined directly on ${\rm SL}(3) \times \R^3 \times \R^3$. The interplay between these
actions is crucial to the analysis; nontrivial isotropy ef\/fectively reduces the number of algebraic relative equilibrium conditions
to be satisf\/ied. Modern treatments of various aspects of this system can be found in~\cite{FassoLewis,Lewis93, ROSD, Rosen98, Rosen88}.
Here we brief\/ly investigate the relative equilibrium
conditions for Lagrangian or Hamiltonian systems sharing the symmetries of the Riemann ellipsoid problem, focusing on
some classes of conf\/iguration-generator pairs identif\/ied by Riemann~\cite{Riemann}.

The historical development af\/f\/ine self-gravitating f\/luid model is an outstanding example of the evolution of a mathematical
model based on qualitative behavior, rather than detailed empirical knowledge. Key qualitative features of the
af\/f\/ine model had been used by Newton and others to predict and interpret the behavior of astronomical bodies long before
 Dirichlet \cite{Dirichlet} rigorously established the af\/f\/ine motions as an invariant subsystem of the PDE modeling the motion
 of an incompressible inviscid f\/luid. The qualitative features of the Riemann ellipsoids has led to their use as models of
 rapidly rotating nuclei \cite{Rosen01, Rosen98, Rosen88}.
Traditionally, selection of an appropriate dynamical model is the f\/irst step in the mathematical analysis of a system;
after the model is specif\/ied, its relative equilibria are identif\/ied and classif\/ied, but
when constructing a model of data for which the key governing phenomena are not known a priori,
but for which certain symmetries are believed to be crucial or certain steady motions have been observed,
we may seek dynamical systems capturing these features. In this situation, it may be desirable to identify
all systems within a diverse family that have the relevant features. Characterizing such systems as elements of
the kernel of an appropriate operator yields computationally tractable conditions that can be
combined with standard data-f\/itting techniques (see Proposition~\ref{rel_crit_conds} and Corollary~\ref{rel_crit_conds_cor}
in Section~\ref{rel_crit_pts}).

A crucial consideration when considering possible applications of modeling of steady motions is the role of
temporally and spatially local approximations. The model need not be a good match to the actual system over extended
periods of time~-- the simplif\/ied model may only be intended to mimic a brief window of behavior of some key subsystem.
For example, a biological or robotic limb may rotate about a joint only through a limited range of angles, and thus
be incapable of steady rotations. However, in assessing the smoothness of motion of the limb,
insight can be gained from measuring the extent to which the actual motion dif\/fers from an interval of
rotation with constant angular velocity about that joint. Assessments of ef\/fectiveness of medication or physical therapy
in improving symptoms of osteoarthritis, or procedures for detection of excessive friction or loading in a robotic arm, can be designed
and implemented even when development of a realistic model of the full system is infeasible.
Conventional statistical analyses can measure the extent to which two motions dif\/fer by comparing measurements
at some key locations, but those measurements may not be well-suited to interpretation and assessment without
further processing. If
desirable motions of the system can be approximated by short term `steady' motions that can be f\/it, e.g.\  using
nonlinear least squares methods, to a family of relative equilibria, then smoothness and correctness of motion can
be quantif\/ied by determination of the goodness of f\/it of measurements to that family.

\section{Relative equilibria as critical points}

The general technique for variational characterizations of relative equilibria can be brief\/ly described as follows:
Consider a manifold $P$ with a $G$ action and a $G$-equivariant vector f\/ield $X$.
Given a point $p_e \in P$ and an element $\bsxi$ of the Lie algebra $\fg$ of
$G$, we wish to determine conditions on $p_e$ and $\bsxi$ under which the
solution of the initial value problem $\dot p = X(p)$ and $p(0) = p_e$ is
\begin{gather}\label{steady_motion}
p(t) = \exp (t   \bsxi) \cdot p_e.
\end{gather}
Such trajectories are known as {\em relative equilibria} or {\em steady motions}. Dif\/ferentiating
(\ref{steady_motion}) with respect to time, we see that a point $p$ is a relative equilibrium with
generator $\bsxi$ of a dynamical systems with vector f\/ield $X$ if\/f $X(p) = \bsxi_P(p)$, where the
vector f\/ield
\begin{gather}\label{inf_gen}
\bsxi_P(p) := \ddt {}\exp (t   \bsxi) \cdot p|_{t = 0}
\end{gather}
is the {\em infinitesimal generator} associated to the algebra element $\bsxi$.
At a  relative equilibrium $p_e$  with continuous isotropy subgroup the generator $\bsxi$ is not unique;
$\bsxi$ may be replaced by $\bsxi + \bszeta$ for any element $\bszeta$ of the isotropy
algebra of $p_e$. Both continuous and discrete isotropy play a~crucial role in the analysis of
relative equilibria (see, e.g., \cite{ Lewis93,Lewis92, ROSD}).

\subsection{Phase space critical point characterizations}

We f\/irst brief\/ly review three well-known characterizations of relative equilibria as critical points of functions
on the full phase space,
beginning with the most traditional and explicit, and moving in the direction of increasing generality.
\begin{itemize}\itemsep=0pt
\item
Energy-momentum method (Lagrangian version).
Consider a Lagrangian system with conf\/iguration manifold $Q$, symmetry group $G$
acting on $Q$, and Lagrangian $L: TQ \to \R$ invariant under the lifted action of $G$ on $TQ$.
The energy of the system is given by the function $E : TQ \to \R$ def\/ined by
\begin{gather*}%\label{energy}
E(v_q) := \FL(v_q)(v_q)  - L(v_q) =  \dep {L((1 + \eps)   v_q) - \eps   L(v_q)} .
\end{gather*}
Here $\FL : TQ \to T^*Q$ denotes the Legendre transformation
\[
\FL(v_q)(w_q) := \dep {L(v_q + \eps   w_q)}.
\]
Noether's theorem gives the momentum map $J: TQ \to \fg^*$ satisfying
\[
J(v_q) \cdot \bsxi = j_\bsxi(v_q) := \FL(v_q)(\bsxi_Q(q))
\]
for all $v_q \in TQ$ and $\bsxi \in \fg$.
The vector f\/ield determined by the Euler--Lagrange equations with Lagrangian $j_\bsxi$ is the inf\/initesimal generator
$\bsxi_{TQ}$ associated to $\bsxi \in \fg$ and the lifted action of $G$ on $TQ$ by~(\ref{inf_gen}).
Thus $v_q$ is a relative equilibrium with generator $\bsxi$ if\/f $v_q$ is a~critical point of $E - j_\bsxi$
(see, e.g., \cite{AbrahamMarsden, MarsdenRatiu}).

\item
Energy-momentum method (symplectic version).
Let $(P, \omega)$ be a symplectic manifold and~$G$ be a Lie group acting on $P$ by symplectomorphisms.
A {\it momentum map} $J: P \to \fg^*$ for the $G$ action is a map such that  for any  $\bsxi \in \fg$
the scalar function $j_\bsxi$ given by $j_\bsxi(p) := J(p) \cdot \bsxi$ has Hamiltonian vector f\/ield equal to the
inf\/initesimal generator $\bsxi_P$ of $\bsxi$.  If there is a momentum map $J$, then $p$ is a relative equilibrium
for a Hamiltonian $H$ if\/f $p$ is a
critical point of $H - j_\bsxi$ (see, e.g., \cite{AbrahamMarsden, GuilStern, LMS, LS89b, LS89a, LMSP89}).

\item
Energy-momentum-Casimir method.
Let $P$ be a Poisson manifold with Poisson bracket $\{\ , \ \}$, and let $G$ act on $P$ by Poisson maps.
A non-constant smooth function $C: P \to \R$ is said to
be a {\it Casimir} if $\{ F,  C \} = 0$ for all smooth functions $F$ on $P$. If there is a
momentum map $J: P \to \fg^*$ associated to the action of $G$ on $(P, \{\ , \ \})$, i.e.\ $\{ F,  j_\bsxi \} = \bsxi_P[F]$
for all smooth functions $F$ on $P$, where $j_\bsxi(p) := J(p) \cdot \bsxi$ as before, then
$p$ is a relative equilibrium with generator $\bsxi$ if\/f $p$ is a critical point of $H - j_\bsxi + C$ for some Casimir $C$
(see, e.g.,  \cite{HHMM, MarsdenRatiu,MarsdenRatiu86, MarsdenWeinstein}).
\end{itemize}

Many Poisson manifolds and noncanonical symplectic manifolds in geometric mechanics arise through the
process of reduction. In symplectic reduction, the reduced phase space associated to the
momentum value $\bsmu \in \fg^*$ is def\/ined as the quotient of the $\bsmu$
level set of the momentum map by the subgroup of group elements preserving
$\bsmu$. The symplectically reduced phase space
has an induced symplectic structure; the reduced equations of motion are
def\/ined using this symplectic structure and the induced Hamiltonian (see, e.g.,
 \cite{AbrahamMarsden, MarsdenWeinstein,Meyer, ZombroHolmes}).
 Poisson reduction proceeds as follows: Given a Poisson manifold $P$ with a free and proper $G$ action by
Poisson maps, the Poisson bracket on the quotient
manifold $M := P/G$ with projection map $\pi: P \to M$ is determined as follows: for any smooth real-valued
functions $f$ and $h$ on $M$, require that
\[
\lcb f, h \rcb_M \circ \pi = \lcb f \circ \pi, h \circ \pi \rcb_P
\]
(see, e.g., \cite{HMRW, MarsdenRatiu}).
 Complications arise in reduction if the group
 action is not free. Continuous isotropy leads to conical singularities of the momentum level sets~\cite{AMM}, and singularities of the quotient manifold arise even in the presence of discrete isotropy.
 Various approaches to addressing these dif\/f\/iculties have been developed: Singular reduction has
 been implemented using stratif\/ied symplectic spaces
 \cite{CushmanSjamaar,OrtegaRatiu, SjamaarLerman} and spaces of invariant functions \cite{ACG, AGJ}.
 In situations where the isotropy is due to the interaction of multiple components of
 a product group action~-- e.g., body and spatial symmetries of a rigid body or continuum system~-- one or
 more components of the group can be excluded from the reduction process and subsequently incorporated
 in the analysis via application of the energy-momentum method to the reduced manifold.

Reduction methods are valuable analytic tools in the analysis of relative equilibria,
but are not always convenient to implement numerically or symbolically. On
a fully symplectically reduced phase space, relative equilibria are critical points of the
Hamiltonian; hence there is no need to match appropriate Casimir
functions to steady motions.  However, there no longer is a single phase space
for the entire system. Rather, there is a distinct reduced phase space
for each value of the momentum and these spaces can change dramatically
as the momentum value is varied. The analysis of these changes can be
quite complicated, and require a substantial amount of topological and
geometric sophistication, even for very low dimensional, simple systems. In the
presence of nontrivial isotropy, these complications are even greater.
The reduced phase space is no longer a manifold, but instead is a stratif\/ied
symplectic space. The variational analysis of relative equilibria is
carried out stratum by stratum. If partial reduction is performed, resulting in a noncanonical structure with
a momentum map associated to the remaining group action, the inequitable treatment of the components
of the symmetry group can obscure the relationship between generators and steady motions, particularly
with points which had nontrivial isotropy in the original manifold.

\subsection{Relative equilibrium characterizations on the conf\/iguration manifold}

Lagrangian systems and their canonical Hamiltonian analogs on (co)tangent bundles with lifted actions have alternative
critical point characterizations. The key observation motivating the reduced energy momentum (REM) method and related approaches
is very simple: due to the structure of the inf\/initesimal generators determined by the  lifted action, the velocity
of a relative equilibrium $v_q \in TQ$ with generator $\bsxi$ satisf\/ies $v_q = \bsxi_Q(q)$; the associated canonical Hamiltonian
systems on $T^*Q$ has relative equilibrium of the form $\FL(\bsxi_Q(q))$, where $\FL$ denote the f\/iber derivative of the Lagrangian
$L$. Thus once an algebra element $\bsxi$ has been selected, there is no need to explicitly solve for the velocity~-- all relative equilibria
with generator $\bsxi$ are uniquely determined by the base point $q$. Given an element $\bsxi$ of the Lie algebra $\fg$ (respectively
 element~$\bsmu$ of~$\fg^*$), we can use a variational  criterion to derive conditions on conf\/igurations guaranteeing that they
are basepoints of relative equilibria with generator~$\bsxi$ (respectively, determine a state with momentum map value~$\bsmu$).

Smale's analysis of relative equilibria of simple mechanical systems with free actions
utilizes a characterization of points in the cotangent bundle generated by inf\/initesimal group motions as
minimizers of kinetic energy within a given level set of the momentum map \cite{SmaleI, SmaleII}. This approach
yields a variational characterization of relative equilibria of such mechanical systems on the conf\/iguration manifold $Q$.
Assume a simple mechanical system, i.e.\ a Lagrangian or canonical Hamiltonian system with conf\/iguration manifold $Q$,
and energy $E(v_q) = \smallfrac 1 2 \norm{v_q}^2 +V(q)$ that is invariant under the lift of a group action on $Q$
\cite{SLM91, SmaleI}, with a free group action.
Given $\bsmu \in \fg^*$, the \textit{amended potential} $V_\bsmu: Q \to \R$ associated to $\bsmu$ is given by
\begin{gather*}%\label{amended_potential}
V_\bsmu(q) = \smallfrac 1 2 \bsmu \cdot \bbI(q)^{-1} \bsmu - V(q),
\end{gather*}
where $\bbI(q): \fg \to \fg^*$ is given by
\[
(\bbI(q) \bsxi) \cdot \bseta := \la \bsxi_Q(q), \bseta_Q(q) \ra_q
\]
and $\bsmu = \bbI(q) \bsxi$. (Since the action is free, $\bbI(q)$ is invertible.)
This explicit construction is linked to Smale's constrained energy minimization characterization as follows:
if we def\/ine $v_\bsmu : Q \to J^{-1}(\bsmu)$ by
\[
v_\bsmu(q) := (\bbI(q)^{-1} \bsmu)_Q(q),
\qquad \mbox{then} \qquad
V_\bsmu= E \circ v_\bsmu.
\]
Critical points $\bq$ of $V_\bsmu$ correspond to relative equilibria $v_\bsmu(\bq)$;
local minima (modulo symmetries of the equilibrium momentum $\bsmu$) of $V_\bsmu$
correspond to (formally) stable relative equilibria \cite{LS89b, SLM91}.

The strategy underlying the REM method for simple mechanical systems with free actions can be extended to more general
functions and symmetries, relaxing the assumptions that the Lagrangian have the form
`kinetic $-$ potential' and that the action be free.
Given $q \in Q$ and $\bsxi \in \fg$, $\bsxi_Q(q)$ is a relative equilibrium of the Euler--Lagrange system
with $G$-invariant regular Lagrangian $L: TQ \to \R$ if\/f $q$ is a critical point of the {\em locked Lagrangian}
$L_\bsxi := L \circ \bsxi_Q$; see \cite{Lewis93, Lewis92}.
One of the key advantages of the (relative) energy-momentum method relative to singular symplectic reduction techniques is the ease
with which continuous isotropy is managed. The locked Lagrangians have no
singularities, even at points with jumps in isotropy; since the momentum level sets do not play an explicit role in
the analysis, possible singularities of the level sets are not a concern.
Analogs of these decompositions for the special case of simple mechanical systems
can be found in \cite{ROlmos}; comparison of the stability conditions obtained in \cite{Lewis93,Lewis92, ROlmos}
illustrates the point that in many situations the geometric setting of steady motion~-- in this case, Lagrangian dynamics on one hand
and cotangent bundle symplectic formulations of simple mechanical systems on the other~-- plays a relatively minor role
in the stability analysis.

\section{Critical points of invariant functions}

Motivated by the considerations touched on in the previous section, we consider critical points of families of
functions parametrized by a Lie algebra or its dual, ignoring any underlying geometric structure beyond a group action.
If the symmetry group is abelian, we can use parametrized families of real-valued group invariant functions.
If the group action is nonabelian, we replace the parametrized family of functions with a single function
depending on both the given manifold and the algebra (or its dual) of the symmetry group, requiring that this function
be invariant with respect to the action on the product manifold determined by the original action and the (co)adjoint action
on the algebra or its dual. We focus on situations in which these functions can be expressed as the composition of an
invariant map to an intermediate manifold with some scalar function on that manifold. This facilitates the identif\/ication of
common properties of a~wide class of functions sharing common symmetry properties.

We begin with two familiar, elementary examples.

\subsection{The Lagrange top}
\label{lag_top}

The heavy Lagrange top is one of the
simplest mechanical systems with a nonlinear conf\/iguration manifold  and a nonfree group action.
Hence we take this system as our f\/irst example of the strategies described above.
The total energy of the heavy top is invariant under spatial rotations about the axis of gravity
and body rotations about the axis of symmetry of the top
(see, e.g.,  \cite{Lewis93, LRSM, Routh} and references therein).

A heavy sleeping  top, i.e.\ a top rotating at constant speed under the inf\/luence of
gravity, with the spatial axis of symmetry aligned with the direction of
the force of gravity, is f\/ixed by the combination of a spatial rotation
through an angle $\theta$ and a body rotation through the angle~$- \theta$. Note
that spin and precession are not uniquely determined for a sleeping
top~-- spin about the axis of symmetry and precession about the axis
of gravity are indistinguishable when those axes are aligned.
This continuous isotropy signif\/icantly complicates many geometric approaches.
Here we brief\/ly review the calculation of the relative
equilibria of the Lagrange top using the locked Lagrangian; see~\cite{Lewis92} for more details.

The Lagrange top has Lagrangian
\begin{gather}\label{lag_top_lag}
L(A   \widehat \bsOmega) = \half \bsOmega^T \bI \bsOmega - g   m   \be_3^T A \bssigma,
\end{gather}
where
\begin{itemize}\itemsep=0pt
\item
$\widehat {\ }: \R^3 \to \{\mbox{$3 \times 3$ skew-symmetric matrices}\}$ is given by $\widehat \bseta   \bszeta = \bseta \times \bszeta$
for all $\bszeta \in \R^3$,
\item
$\bI = I_1   \idm + (I_3 - I_1) \bssigma \bssigma^T$ is the inertia tensor of the top, for positive constants
$I_1$ and $I_3$, and unit vector $\bssigma$ determining the axis of symmetry of the top,
\item
$g  \be_3$ is the force of gravity,
\item
$m   \bssigma$ is the center of mass.
\end{itemize}
Note that if we identify the Lie algebra ${\rm so}(3)$ with $\R^3$ using $\widehat {\ }$\,, then
the left trivialization $T {\rm SO}(3) \to {\rm SO}(3) \times \R^3$ of the tangent bundle takes the form
$A   \widehat \bsOmega \mapsto (A, \bsOmega)$.

The Lagrangian (\ref{lag_top_lag}) is invariant under the action of $S^1 \times S^1$  on ${\rm SO}(3)$ by
\begin{gather}\label{top_action}
(\theta, \varphi) \cdot A = \exp(\theta   \be_3) A \exp(-\varphi   \bssigma),
\end{gather}
corresponding to spatial rotations about the vertical axis $\be_3$ and body rotations about the axis of symmetry
$\bssigma$.

The locked Lagrangian $L_{(\xi_\ell, \xi_r)}: {\rm SO}(3) \to \R$ determined by $L$ and the algebra element
$(\xi_\ell, \xi_r) \in \bR^2 \approx \fg$ has the form
\[
L_{(\xi_\ell, \xi_r)}(A) = L\big(A(\xi_\ell   \widehat \be_3 - \xi_r   \widehat A \bssigma)\big)
 = \half (\xi_\ell   \be_3 - \xi_r   A \bssigma)^T \bI (\xi_\ell   \be_3 - \xi_r   A \bssigma) - g   m   \be_3^T A \bssigma.
 \]
$L_{(\xi_\ell, \xi_r)}$ is a quadratic polynomial in the invariant $\iota: {\rm SO}(3) \to [-1, 1]$ given by
\begin{gather}\label{heavy_top_inv}
\iota(A) := \be_3^T A \bssigma,
\end{gather}
with dif\/ferential
\[
d\iota(A)(A \widehat \bseta) = \be_3^T \widehat \bseta A \bssigma = (\be_3 \times A \bssigma)^T \bseta.
\]
Thus the critical points of $\iota$ are the `sleeping' states, with $A \bssigma = \pm \be_3$; the
critical point set is $\iota^{-1}(\{-1, 1\})$, the preimage of the endpoints of the codomain of $\iota$.
Sleeping states are f\/ixed by the $S^1$ action. All other points have trivial isotropy.

The locked Lagrangian satisf\/ies $L_{(\xi_\ell, \xi_r)} = \ell_{(\xi_\ell, \xi_r)} \circ \iota$, where
\[
\ell_{(\xi_\ell, \xi_r)}(x) := \half \xi_\ell^2 (I_3 - I_1) x^2 - (\xi_\ell   \xi_r   I_3 + g   m) x
+ \half \big(\xi_\ell^2 I_1 + \xi_r^2 I_3\big).
\]
Clearly, sleeping tops are critical points of $L_{(\xi_\ell, \xi_r)}$ for any parameter values.
`Tilted' states, i.e.\ elements of $\iota^{-1}((-1, 1))$, are critical points of $L_{(\xi_\ell, \xi_r)}$ if\/f
\begin{gather}\label{tilt_eq}
\xi_\ell^2 (I_3 - I_1) \iota(A) = \xi_\ell   \xi_r   I_3 + g   m.
\end{gather}
For any $A \in {\rm SO}(3)$ and  nonzero value of $\xi_\ell$, there is a unique value of $\xi_r$ such that
$A$ is a critical point of $L_{(\xi_\ell, \xi_r)}$, and hence determines a relative equilibrium with generator $(\xi_\ell, \xi_r)$.

The only true equilibria are the sleeping states. Note that $\ell_{(0, 0)}(x) = - g   m  x$ is linear; hence the critical points of
$L_{(0,0)}$ are precisely those of $\iota$. The Lagrangian $L$ is invariant with respect to negation of the velocity, corresponding
to the time reversal symmetry of the dynamics; this invariance induces the `symmetry' that two algebra elements dif\/fering
only by sign determine the same locked Lagrangian. The origin is the unique f\/ixed point of the $\bZ_2$ action on $\R^2$,
$\bsxi \mapsto - \bsxi$. Of course, in the case at hand, it is hardly surprising that the conditions for true equilibrium are signif\/icantly
 simpler than those for relative equilibrium; however, we shall see that in systems with more elaborate symmetry groups, consideration
 of the relevant group action on the algebra and classif\/ication of the algebra elements with regard to isotropy can be a powerful
 tool in the analysis of relative equilibria.

\begin{remark}
An upright state, with $\iota(A) = 1$, has isotropy algebra $\mbox{span}\{(1, 1)\}$;
setting $(\xi_\ell, \xi_r) = (\lambda, \lambda - \xi)$ in \eqref{tilt_eq} sheds light on the relationship between the
upright relative equilibrium and nearby tilted states with trivial isotropy. Solving for $\lambda$ yields
\[
\lambda = \frac {I_3 \xi \pm \sqrt{\xi^2 I_3^2  - 4 \, g \, m (I_3 + (I_1 - I_3) \iota(A) }} {2 (I_3 + (I_1 - I_3) \iota(A))}.
\]
Thus, setting $\iota(A)$, we see that the family of solutions \eqref{tilt_eq} intersects the family of
sleeping tops throughout the range of angular velocities $\xi^2 I_3^2 >
4 g m I_1$ for which the sleeping tops are stable. Corresponding to every speed for which
the sleeping top rotates stably, there is a family of tilted, precessing
motions whose combined spin and precession equal the angular velocity of the sleeping top \cite{Lewis93}.
Some approaches, particularly those involving stratif\/ied spaces or reduction with respect to a subgroup, can
obscure this simple but signif\/icant relationship between the generators of sleeping and tilted states.
\end{remark}

The division of critical points into sleeping and tilted states holds for any parametrized function of the form
$F = f \circ \iota$ for some function $f: \R \to \R$. Sleeping states are critical points for any
function $f$; additional critical points of $F$ are states $A$ such that $\iota(A)$ is a critical point of~$f$.

\subsection{The double spherical pendulum}
\label{pendulum}

Our second mechanical example with an abelian symmetry group is the double spherical pendulum. In this
example, discrete isotropy subgroups play a central role in the classif\/ication of the relative equilibria.

Let $\bq_j \in S^2$ denote the normalized position of the $j$-th
pendulum, with mass $m_j$ and length~$\ell_j$.
Let $g$ denote the strength of
gravity and $\xi \in \bR$ denote the generator, i.e.\ the (vertical) angular
velocity of the pendula. The Lagrangian for the double spherical pendulum is
\[
L(\dot \bq_1, \dot \bq_2) = \smallfrac {m_1} 2 \norm{\ell_1 \dot \bq_1}^2 + \smallfrac {m_2} 2 \norm{\ell_1 \dot \bq_1 + \ell_2 \dot \bq_2}^2
	- g   \be_3^T( (m_1 + m_2)  \ell_1 \bq_1 + m_2 \ell_2 \bq_2 ).
\]
The relative equilibria of the double spherical pendulum are studied using a combination of Lagrangian
reduction and blow-up techniques. Here, as in our previous example, we use the locked Lagrangian in our analysis.
The locked Lagrangian
\begin{gather}
L_\xi(\bq_1, \bq_2)  =   \smallfrac {\xi^2} 2 \big( m_1 \ell_1^2 \norm{\be_3 \times \bq_1}^2 + m_2 \norm{\be_3 \times (\ell_1 \bq_1 + \ell_2 \bq_2)}^2 \big)
\nonumber\\
\hphantom{L_\xi(\bq_1, \bq_2)  =}{}
	- g  \be_3^T( (m_1 + m_2)  \ell_1 \bq_1 + m_2 \ell_2 \bq_2 ).\label{dsp_ll}
\end{gather}
for the double spherical pendulum is invariant under the diagonal action of
${\rm O}(2)$ consisting of rotations about the vertical axis and ref\/lections across vertical planes. We f\/irst formulate
critical point conditions for a broad class of functions on $S^2 \times S^2$ invariant under the $S^1$ action
described above, then specialize those conditions to the locked Lagrangian for the double spherical pendulum.
We def\/ine the invariant map $\iota: S^2 \times S^2 \to \R^3$ given by
\begin{gather}\label{dsp_inv}
\iota(\bq) := (\be_3 \cdot \bq_1, \be_3 \cdot \bq_2, \bq_1 \cdot \bP_3 \bq_2),
\end{gather}
where $\bP_3$ denotes projection into the horizontal plane, and consider invariant functions of the form
$F = f \circ \iota$ for some function $f: \R^3 \to \R$. Note that if $\{\bq_1, \bq_2, \be_3 \}$ is linearly dependent, then
the components of $\iota = \iota(\bq_1, \bq_2)$ satisfy the relation
\begin{gather}\label{dependence}
\iota_3^2 = \lp 1 - \iota_1^2 \rp \lp 1 - \iota_2^2 \rp.
\end{gather}

The critical point behavior is determined by the isotropy of $\bq = (\bq_1, \bq_2)$, which is in turn determined by  the dimension of
$\mbox{span} \{\bq_1, \bq_2, \be_3 \}$:
\begin{proposition}\qquad
\label{pend_prop}

\begin{itemize}\itemsep=0pt
\item
If $\bq_1$ and $\bq_2$ are both vertical, then $\bq$ has full isotropy and $\bq$ is a relative critical point of any invariant function.

\item
If $\{\bq_1, \bq_2, \be_3 \}$ spans a plane, then $\bq$ is fixed by reflection across that plane.

If $|\iota_j| = 1 \neq |\iota_{j'}|$, where $\{j, j'\}$ is  a permutation of $\{1, 2\}$, $\bq$ is a critical point iff
\[
\pd f {\iota_{j'}}(\iota) =   \pd f {\iota_3}(\iota) = 0.
\]
If $|\iota_j| \neq 1 \neq |\iota_{j'}|$, $\bq$ is a critical point iff
\begin{gather}\label{plane_span}
\lp 1 - \iota_j^2 \rp \pd f {\iota_j}(\iota)
	 = \iota_j \iota_3 \pd f{\iota_3}(\iota), \qquad   \mbox{for} \quad j = 1, 2 \quad \mbox{and} \quad \iota = \iota(\bq).
\end{gather}

\item
If $\{\bq_1, \bq_2, \be_3 \}$ is linearly independent, then $\bq$ has trivial isotropy and
is a critical point iff~$\iota(\bq)$ is a critical point of $f$.
\end{itemize}
\end{proposition}

\begin{proof}
If we def\/ine $\bsrho_j := \bq_j \times \be_3$, $j = 1, 2$,  then
\begin{gather*}%\label{diota_tsp}
d \iota(\bq)=   \begin{pmatrix}
\bsrho_1^T & 0 \\
0 & \bsrho_2^T \\
(\bq_1 \times \bP_3 \bq_2)^T & (\bq_2 \times \bP_3 \bq_1)^T
\end{pmatrix}.
\end{gather*}

If $\bq_j$ is vertical, then $|\iota_j| = 1$, $\iota_3 = 0$, and $\bsrho_j = 0$. If both $\bq_1$ and $\bq_2$ are vertical, then
$d \iota(\bq) = 0$; hence $\bq$ is a critical point of any invariant function. If $\bq_1$ is vertical, but $\bq_2$ is not,
then
\[
d \iota(\bq)=  \begin{pmatrix}
0 & 0 \\
0 & \bsrho_2^T \\
-\iota_1 \bsrho_2^T & 0
\end{pmatrix}.
\]
Since $\bsrho_2 \in T_{\bq_j} S^2$ for $j = 1, 2$ in this case, $\mbox{range} [d \iota(\bq)] = \mbox{span}\{\be_2, \be_3 \}$.
The case in which $|\iota_2| = 1 \neq |\iota_1|$ is analogous.

If neither pendulum is vertical, we can take
\[
\lcb   \begin{pmatrix} \bsrho_1 \\ 0 \end{pmatrix},   \begin{pmatrix} 0 \\ \bsrho_2 \end{pmatrix}  ,
  \begin{pmatrix} \bsrho_1 \times \bq_1\\ 0 \end{pmatrix}  ,   \begin{pmatrix} 0 \\ \bsrho_2 \times \bq_2  \end{pmatrix}  \rcb
\]
as a basis for the tangent space. With respect to this basis, $d \iota(\bq)$ has matrix representation
\begin{gather}\label{mat_rep}
  \begin{pmatrix}
\norm{\bsrho_1}^2 & 0 & 0 & 0 \\
0 & \norm{\bsrho_2}^2 & 0 & 0 \\
- \iota_1 \iota_3 & - \iota_2 \iota_3 & \tau & - \tau
\end{pmatrix}  ,
\qquad   \mbox{where} \qquad \tau = \be_3^T (\bq_1 \times \bq_2).
\end{gather}
Hence $d \iota(\bq)$ has full rank if\/f $\tau \neq 0$, i.e.\ $\{ \bq_1, \bq_2, \be_3 \}$ is linearly independent.
(\ref{plane_span}) now follows directly from consideration of the nullspace of the transpose of (\ref{mat_rep}).
\end{proof}

We now use Proposition \ref{pend_prop} to f\/ind the relative equilibria of the double pendulum system.
We assume that the
 pendulum shaft lengths $\ell_j$ and masses $m_j$ are nonzero for $j = 1, 2$, and that $g \neq 0$.
 Note that the locked Lagrangian (\ref{dsp_ll}) for the double spherical pendulum  is
 smooth on all of $S^2 \times S^2$, while the version of the amended potential used in \cite{MarsdenScheurle}
 is not def\/ined on the vertical conf\/igurations and not dif\/ferentiable at states for which at least one of the
 pendula is horizontal.  The calculations given below recapture the results of~\cite{MarsdenScheurle}, and additional
 relative equilibria in which one of the pendula is horizontal; relative equilibrium with both pendula horizontal
 is possible only if the shaft lengths are equal.

\begin{proposition}\label{dsp_prop}
If $\bq = (\bq_1, \bq_2)$ and $\xi$ determine a relative equilibrium of the double spherical pendulum, then either both
or neither of $\bq_1$ and $\bq_2$ are vertical.

The four configurations in which both pendula are vertical are the only equilibria.

If neither  $\bq_1$ nor $\bq_2$ is vertical, then $\bq = (\bq_1, \bq_2)$ and $\xi \neq 0$ determine a relative equilibrium iff
$\bq_1$ and $\bq_2$ lie in a common vertical plane
and the parameter ratios
\[
\ell := \bigfrac {\ell_2}{\ell_1}, \qquad m := \bigfrac {m_2 \, \ell}{m_1 + m_2}, \qquad {\rm and} \qquad  \gamma := \frac g {\xi^2 \ell_1}
\]
and components of $\iota = \iota(\bq)$ satisfy one of the following sets of conditions
\begin{itemize}\itemsep=0pt
\item
$\displaystyle \iota_1 = \frac {1 - \gamma}{1 +  \frac m {\varkappa}}$
and
$\displaystyle \iota_2 = \frac {\ell - \gamma}{\ell  + \varkappa}$,
 where  ${ \varkappa := \frac {\iota_3} {1 - \iota_2^2}}$,   or

\item
$\gamma = r_j$ (equivalently, $\ell_j \xi^2 =  g$),  ${\displaystyle \iota_{j'}  = \frac  {1 - s_j}{1 - \frac m \ell}}$,
${\displaystyle \iota_j^2 = 1 - t_j^2 \lp  1- \iota_{j'}^2 \rp}$, and $\iota_3 < 0$, \\
for $\{j, j'\} = \{1, 2\}$ or $\{2, 1 \}$, where $\br = (1, \ell)$, $\bs = \lp \smallfrac 1 \ell, \ell \rp$, and $\bt = \lp m, \smallfrac 1 \ell \rp$.
\end{itemize}
\end{proposition}

Note that the parameter $\ell$ corresponds to the ratio $r$ in~\cite{MarsdenScheurle}, while the parameter $m$ used here is the product of
$\ell$ and the reciprocal of the parameter also denoted $m$ in~\cite{MarsdenScheurle}.

\begin{proof}
The locked Lagrangian (\ref{dsp_ll}) for the trivial generator is linear in the invariant:
\[
L_0(\bq) = - g (m_1 \ell_1, m_2 \ell_2, 0)^T \iota(\bq).
\]
The only critical points of $L_0$ are the critical points of $\iota$, namely the vertical pendula.

For $\xi \neq 0$, the locked Lagrangian $L_\xi$  satisf\/ies $L_\xi = \xi^2 \ell_1^2 (m_1 + m_2) f \circ \iota$ for
\begin{gather}\label{pend_f}
f(\iota) := \half  \lp 1 - \iota_1^2 + m   \ell \lp 1 - \iota_2^2 \rp + m   \iota_3 \rp - \gamma (\iota_1 + m   \iota_2),
\end{gather}
The function (\ref{pend_f})  has dif\/ferential
\[
d f (\iota) =   \begin{pmatrix}
1 - \iota_1 - \gamma \\
m \lp \ell (1 - \iota_2) - \gamma \rp \\
m
\end{pmatrix}
\]
Thus, since $\pd f {\iota_3} \equiv m \neq 0$, Proposition \ref{pend_prop} implies that if $\bq$ is a  critical point
of the locked Lagrangian, then $\bq_1$ and $\bq_2$ lie in a common vertical plane;
either both pendula are vertical or neither is vertical. In the case that neither pendulum is vertical and $\{ \bq_1, \bq_2, \be_3 \}$
is dependent, the criticality conditions (\ref{plane_span}) take the form
\begin{gather}
 1 - \iota_1 - \gamma = \pd f {\iota_1}(\iota) = \frac {\iota_1} \varkappa  \pd f {\iota_3}(\iota) = \frac {m   \iota_1} \varkappa,\nonumber \\
  m (\ell ( 1 - \iota_2) - \gamma) = \pd f {\iota_2}(\iota) = \iota_2 \varkappa    \pd f {\iota_3}(\iota) = m   \iota_2  \varkappa. \label{pend_eq_cons}
\end{gather}

We can decompose the analysis of (\ref{pend_eq_cons}) into subcases determined by the relationships between the parameters $\ell$ and $\gamma$:

$\bullet$ $1 \neq \gamma \neq \ell$.
The criticality conditions (\ref{pend_eq_cons}) have the unique solution
\begin{gather}\label{non_special_sp}
\iota_1 = \frac {1 - \gamma}{1 + \frac m \varkappa}
\qquad \text{and}\qquad
\iota_2 = \frac {\ell - \gamma}{\ell + \varkappa}.
\end{gather}

$\bullet$
$\gamma  = 1$.
The f\/irst equation in (\ref{pend_eq_cons}) simplif\/ies to $\lp 1 + \frac m \varkappa \rp \iota_1 = 0$, which is satisf\/ied if\/f $\iota_1 = 0$ or
$-\varkappa = m$. The solution $\iota_1 = 0$ corresponds to the solution (\ref{non_special_sp}).
If $- \varkappa =  m$, then the f\/irst equation in (\ref{pend_eq_cons}) is satisf\/ied for any value of $\iota_1$, and the second equation can
be rearranged as $\iota_2 = \frac {\ell - \gamma}{\ell - m}$. (Note that $m_1 > 0$ implies that $\ell > m$.)
\[
\iota_3 = \varkappa \lp 1 - \iota_2^2 \rp = - m \lp 1 - \iota_2^2 \rp < 0.
\]
Finally, $\iota_1$ is determined up to sign by the equality
\begin{gather}\label{varkappa}
m^2 = \varkappa^2 = \frac {1 - \iota_1^2}{1 - \iota_2^2}
\end{gather}
if $\bq_1$ and $\bq_2$ lie in a common vertical plane and the hypothesis $\varkappa = - m$ is satisf\/ied.

$\bullet$
$\ell =  \gamma$.
The second equation in (\ref{pend_eq_cons}) simplif\/ies to $\lp \ell + \varkappa \rp \iota_2 = 0$, which is satisf\/ied if\/f $\iota_2 = 0$ or
$-\varkappa = \ell$. The derivation of the conditions on $\iota$ is analogous to that for the case $\gamma = 1$.
\end{proof}

To identify orbits of relative equilibria with generator $\xi$ and image $\iota$ under the invariant map, we can use the relationship
(\ref{varkappa}) to formulate the search for relative equilibria $(\bq, \xi)$ satisfying the f\/irst set of conditions in Proposition~\ref{dsp_prop}
as a search for f\/ixed points of a pair of maps.
If we introduce the mappings $\kappa: (-1, 1)^2 \to \R^+$ and $\Upsilon_\pm: (-1, 1)^2 \to \R^2$ given by
\[
\kappa(\iota_1, \iota_2) := \sqrt{\frac {1 - \iota_1^2}{1 - \iota_2^2}}
\qquad {\rm and}\qquad
\Upsilon_\pm(\iota_1, \iota_2) := \lp \frac {1 - \gamma}{1 \pm \frac m { \kappa(\iota_1, \iota_2)}},  \frac {\ell - \gamma}{\ell \pm  \kappa(\iota_1, \iota_2)} \rp,
\]
then $\iota$ satisf\/ies (\ref{non_special_sp}) for $1 \neq \gamma \neq \ell$ if\/f $(\iota_1, \iota_2)$ is a f\/ixed point of $\Upsilon_{{\rm sgn}(\iota_3)}$.
The following is one possible approach to f\/inding these f\/ixed points: Substituting the expressions (\ref{non_special_sp}) into (\ref{varkappa}) and
regrouping terms yields a sixth order polynomial in $\varkappa$ whose
roots determine candidate values of~$\iota_1$ and~$\iota_2$; only those values for which $|\iota_j| \leq 1$, $j = 1, 2$ are
admissible. Note that this polynomial can be regarded as a quadratic polynomial in~$\gamma$, with~$\varkappa$ treated as a
parameter if it is preferable to solve for the angular velocity~$\xi$ as a function of~$\iota_1$ and~$\iota_2$.

\subsection{Adequate invariant maps}
\label{adequate}

Analysing classes of invariant functions, rather than a single function or parametrized
family of functions, can shed light on the role of symmetry in determining the critical
point conditions (see, e.g., \cite{Field, GSS, Olver95, Olver94, Olver09}).
This can be particularly valuable in situations where the appropriate
parameterized models are not known a priori and choices must be made on the basis
of a~range of desired behaviors. We outline a possible approach that can be used
to symbolically or numerically search for functions having specif\/ied critical point and
isotropy relationships. Generalizing the approach of the previous section, we consider
invariant functions constructed as the composition of an invariant map and
scalar  functions on the codomain of  that map. If $\iota: P \to \cQ$ is a $G$-invariant map
from $P$ to a manifold $\cQ$ (possibly with boundary) and
$F = f \circ \iota$ for some $f: \cQ \to \R$, then $p$ is a critical point of $F$ if\/f
$d f(\iota(p)) \in \ker \lsb d_p^*  \iota \rsb$.
By characterizing these kernels, we can obtain critical point conditions for such
invariant functions.

Generating sets of invariants are known for many important group actions, and can be used to construct an
invariant map $\iota: P \to \cQ$ (see, e.g., \cite{GSS, KMO, KoganOlver}.)
 If a generating set is not known a priori and cannot be readily determined, more
 limited criteria for diversity of an invariant map $\iota: P \to \cQ$ can be used. Such a map may not capture
 the full range of possible invariant function behavior, but can still be used to study very general classes
 of invariant functions. We describe a possible criterion~-- that the rank of the dif\/ferential of the map be as
 large as possible at all points~-- that is relatively easy to verify in many cases, and then
 illustrate implementation of this test in several examples.

Let $\Phi_g: P  \to P$ denote the action of $g \in G$ on a manifold $P$, i.e.\ $\Phi_g(p) = g \cdot p$.
If $F: P \to \R$ is $G$-invariant, then for any $\bsxi \in \fg$,
\[
\bsxi_P[F] = \ddt {} F \circ \Phi_{\exp (t \, \bsxi)}|_{t = 0}  = \ddt {F}= 0;
\]
hence for any $p \in P$,
$\rang {d_p^* F} \subseteq ( \fg \cdot p)^A$, the annihilator of the tangent space to the group orbit through $p$.
There are additional restrictions on the dif\/ferentials of invariant functions at points with nontrivial isotropy.

\begin{proposition}
If  $\iota: P \to \cQ$  is a differentiable $G$-invariant map $\iota: P \to \cQ$, then
\begin{gather*}%\label{containment}
\rang {d_p^*  \iota} \subseteq  \lp \fg \cdot p \cup \lp  \cup_{g \in G_p} \rang{d_p \Phi_g - \idm} \rp \rp^A,
 \end{gather*}
 where $G_p$ denotes the isotropy subgroup of $p$.
\end{proposition}

\begin{proof}
If $g \in G_p$, the isotropy subgroup of $p$, and $F: P \to \R$ is $G$-invariant, then
\[
dF(p) \circ d_p \Phi_g =  dF(\Phi_g(p)) \circ d_p \Phi_g = d(F \circ \Phi_g)(p)  = dF(p)
\]
implies
\begin{gather}\label{rang_act}
\rang{d_p \Phi_g - \idm} \subseteq \ker [dF(p)],
\qquad \mbox{i.e.} \qquad
dF(p) \in \ker \lsb \lp d_p \Phi_g \rp^* - \idm \rsb.
\end{gather}
If $F = f \circ \iota$, then $d F(p)(\delta p) = df(\iota(p))(d_p \iota(\delta p)) = d_p^*  \iota df(\iota(p))(\delta p)$. Hence
$G$ invariance of $F$ implies that $d_p^*  \iota$ maps dif\/ferentials at $\iota(p)$ of functions on
$\cQ$ into the annihilator of the tangent space $\fg \cdot p$ of the group orbit through $p$.
(\ref{rang_act}) implies that $d_p^*  \iota (df(\iota(p))) \in \ker \lsb \lp d_p \Phi_g \rp^* - \idm \rsb$
if $g \in G_p$, and hence
\[
d_p^*  \iota (df(\iota(p)))  \in
( \fg \cdot p)^A \cap \lp  \cap_{g \in G_p} \ker \lsb d_p^* \Phi_g - \idm \rsb \rp
 =  \lp \fg \cdot p \cup \lp  \cup_{g \in G_p} \rang{d_p \Phi_g - \idm} \rp \rp^A
 \]
 for any dif\/ferentiable function $f: \cQ \to \R$.
\end{proof}

 \begin{definition}
A $G$-invariant map $\iota: P \to \cQ$ is {\em adequate} if
 \begin{gather}\label{adequate_cond}
\mbox{rank} \lsb d_p^*  \iota \rsb
+ \dim  \lp \fg \cdot p \cup \lp  \cup_{g \in G_p} \rang{d_p \Phi_g - \idm} \rp \rp
= \dim P
\end{gather}
for all $p \in P$.
\end{definition}

The following proposition illustrates the process of verifying adequacy in some simple and familiar settings.

\begin{proposition}\label{adequate_examples}\qquad
\begin{enumerate}\itemsep=0pt
\item[$1.$]
$\iota: \R^3 \to \R$ given by $\iota(\bx) := \norm{\bx}^2$ is adequate for the usual action of ${\rm O}(3)$ on $\R^3$.

\item[$2.$]
$\iota: {\rm SO}(3) \to S^2$ given by $\iota(A) := A^T \bssigma$ is adequate for the action of  $S^1$ on ${\rm SO}(3)$ by
$\theta \cdot A = R_\theta A$, where  $R_\theta$ implements rotation through $\theta$ about  a fixed axis $\bssigma \in S^2$.

\item[$3.$]
The map \eqref{heavy_top_inv} is adequate for the action of $S^1 \times S^1$ on ${\rm SO}(3)$ by \eqref{top_action}.

\item[$4.$]
The map \eqref{dsp_inv} is adequate for the action of ${\rm O}(2)$ on $S^2 \times S^2$ given in Section~{\rm \ref{pendulum}}.

\item[$5.$]
$\iota: {\rm GL}(3) \to (\R^+)^3$ given by
\[
\iota(A) :=  \lp \mbox{\rm trace} \, C, \half \lp (\mbox{\rm trace} \, C)^2 - \mbox{\rm trace} \, C^2 \rp\!, \det C \rp,
\qquad   \mbox{where} \quad C := A^T A,
\]
is adequate for the action of ${\rm O}(3) \times {\rm O}(3)$ acting on ${\rm GL}(3)$ by $(U, V) \cdot A = U A V^T$.
Equivalently, setting $\iota_0 = 1$,
\[
\det\big(A^T A  - \lambda   \idm \big) =  \sum_{j = 0}^3 \iota_j(A) (- \lambda)^{3 - j}
\qquad  \forall \, \lambda \in \R.
\]
\end{enumerate}
\end{proposition}
\begin{proof}
1.~$d_\bx^* \iota$ is rank one and $\fg \cdot \bx$
is two-dimensional if $\bx \neq 0$, while $\rang{d_0 \Phi_{- \idm} - \idm} = \rang{- 2 \, \idm} = \R^3$ implies
\[
\cup_{A \in G_{\boldsymbol 0} = {\rm O}(3)} \rang{d_0 \Phi_A - \idm} = \R^3.
\]
Hence $\iota$ is adequate.

2.~$\mbox{range}[d_A\iota] = T_{\iota(A)} S^2$ for all $A \in {\rm SO}(3)$ and $\fg \cdot A$ is one-dimensional for all $A$,
so (\ref{adequate_cond}) is satisf\/ied everywhere and $\iota$ is adequate.

3.~If $A$ is not a critical point of $\iota$, then $d \iota(A)$ is rank one and $A$ has trivial isotropy. It follows that
$\dim \fg \cdot A = \dim \lp S^1 \times S^1 \rp$. Hence $\mbox{rank} [d_A \iota] + \dim \fg \cdot A
= 3$ at these points.

If $A$ is a critical point of $\iota$, and hence $A \bssigma = \pm \be_3$, then
\[
d_A \Phi_{(\theta, \pm \theta)}(\widehat \bseta A) = R_{\theta   \be_3} \widehat \bseta A R_{\mp \theta   \bssigma}
 = R_{\theta  \be_3} \widehat \bseta R_{\mp \theta   A \bssigma} A
  = R_{\theta   \be_3} \widehat \bseta R_{- \theta   \be_3 } A
    = \widehat {R_{\theta   \be_3} \bseta }A
     \]
 implies that
 \[
 \mbox{range}\lsb d_A \Phi_{(\theta, \pm \theta)} - \idm \rsb = \mbox{span} \lcb \widehat \be_1 A,   \widehat \be_2 A \rcb
 \qquad    \mbox{if} \quad \theta \neq 2   \pi   k, \quad k \in \bbZ.
 \]
 Combining this with
 $\fg \cdot A =  \mbox{span} \lcb \widehat \be_3 A, A \widehat \bssigma \rcb = \mbox{span} \lcb \widehat \be_3 A \rcb$
 if $A \bssigma = \pm \be_3$, we see that in this case
 \[
  (\fg \cdot A) \cup \lp \cup_{\theta \in S^1}  \mbox{range}\lsb d_A \Phi_{(\theta, \pm \theta)} - \idm \rsb \rp = T_A {\rm SO}(3).
\]
Hence $\iota$ is adequate.

4.~If we regard $T_{\bq_j} S^2$ as a subspace of $\R^3$, then
\[
d_\bq \Phi_\pi - \idm = - 2   \bP_3 \times \bP_3,
\]
where $\bP_3$ denotes orthogonal projection onto the $\be_1$--$\be_2$ plane. Thus if both $\bq_1$ and $\bq_2$ are
vertical, then $\rang{d_\bq \Phi_\pi - \idm} = T_\bq \lp S^2 \times S^2 \rp$, and hence (\ref{adequate_cond}) is satisf\/ied.
If $\bq_1$ and $\bq_2$ lie in the same
vertical plane, then Proposition~\ref{pend_prop} implies that $\mbox{rank}[d_\bq^* \iota] = 2$; since
\[
\rang{d_\bq \Phi_\pi - \idm} \supseteq
\mbox{span} \{(\be_3 \times \bq_1, \bzero), (\bzero, \be_3 \times \bq_2)\}
\]
in this case, (\ref{adequate_cond}) is satisf\/ied at these points.
If $\mbox{span} \{\bq_1, \bq_2, \be_3 \} = \R^3$, then $\mbox{rank}[d_\bq^* \iota] = 3$ and
$\fg \cdot \bq$ is one-dimensional. Hence (\ref{adequate_cond}) is satisf\/ied at all points and $\iota$ is adequate.

5.~Noting that $\mbox{rank} [d_A \iota] = \#$  distinct singular values of $A$, we now consider case by case matrices
with dif\/ferent numbers of distinct singular values.

If $A$ has three distinct singular values, then the isotropy subgroup of $A$ is discrete,
$\fg \cdot A$ is six-dimensional, and the dif\/ferentials of the components of $\iota$ are linearly independent.

If $A$ has only two distinct singular values, then there are matrices $U, V \in {\rm O}(3)$ and
sca\-lars~$\alpha$,~$\tilde \alpha$ such that
\[
A = U \lp \alpha   \idm + \tilde \alpha \be_3 \be_3^T \rp V^T.
\]
If we set $\bu := U \be_3$ and $\boldv := \det (U \, V)  V \be_3$, then the isotropy subgroup $G_A$ of $A$ is
\[
G_A = \lcb \lp R_{\theta  \bu}, R_{\theta  \boldv} \rp : \theta \in S^1 \rcb,
\]
where $R_\bx$ denotes rotation about $\bx$ through the angle $|\bx|$ for any nonzero $\bx \in \R^3$.
If we def\/ine $\psi_\theta: \R^{3 \times 3} \to  \R^{3 \times 3}$ by
\[
\psi_\theta(B) := R_{\theta   \bu} B R_{\theta  \boldv}^T - B,
\]
then $(\fg \cdot A) \cup \lp \cup_{\theta \in S^1} \mbox{range}\lsb \psi_\theta \rsb \rp $ has orthogonal complement
$\mbox{span}\{U V^T, \bu \boldv^T \}$ with respect to the inner product
$\langle B, \widetilde B \rangle = \mbox{trace}(B^T  \widetilde B)$.
If $T_A^* {\rm GL}(3)$ is identif\/ied with $\R^{3 \times 3}$ by means of this inner product, then
\[
\rang{d_A^* \iota} = \mbox{span}\big\{U V^T, \bu \boldv^T \big\}
\]
as well. Hence (\ref{adequate_cond}) holds for $A$ with multiple  singular values.

If $A$ has exactly one singular value $\sigma$, then $\frac 1 \sigma A \in {\rm O}(3)$ and
\[
G_A = \{ (A U, U\!  A): U \in {\rm O}(3) \}.
\]
The linearization of the action of an element of the isotropy subgroup satisf\/ies
\[
d_A \Phi_{(A U, U   A)} (V   A) = A U ( V   A )(U  A)^T  = A U V (A U)^T A = ( \Ad {A U} V) A.
\]
Hence if $\rho_A$ denotes right multiplication by $A$, then
\[
d_A \rho_A^{-1} \circ \lp d_A \Phi_g - \idm_{T_A {\rm GL}(3)} \rp \circ d_\idm \rho_A
= \Ad {A U} - \idm_{\R^{3 \times 3}}
\]
and thus
\begin{gather*}
d_A \rho_A^{-1} \lp \cup_{g \in G_A} \rang{d_A \Phi_g - \idm_{T_A {\rm GL}(3)}} \rp\\
\qquad{}
 =   \cup_{U \in {\rm O}(3)} \rang{\Ad {A U} - \idm_{\R^{3 \times 3}}}
 =  \{\mbox{traceless $3 \times 3$ real matrices} \}.
\end{gather*}
If $\frac 1 \sigma A \in {\rm O}(3)$, then $\rang{d_A^* \iota} = \mbox{span}\{\mbox{trace} \circ d_A \rho_A^{-1} \}$. Hence
(\ref{adequate_cond}) is satisf\/ied for $A$ with a~unique singular value and $\iota$ is adequate.
\end{proof}

The f\/irst three examples in Proposition~\ref{adequate_examples} are well-known quotient maps;
the second example illustrates the possible advantages of regarding $\iota$ as a map into a manifold, rather
than a collection of invariant scalar functions. The fourth example, the ${\rm O}(2)$ action on $S^2 \times S^2$, illustrates
the distinction between an invariant functions approach and a quotient space approach: the components
of the map (\ref{dsp_inv}) are related by (\ref{dependence}) when $\{\bq_1, \bq_2, \be_3 \}$ is linearly dependent,
as they must be if either the f\/irst or second component of $\iota(\bq)$ has unit absolute value.
(As we discuss brief\/ly below, such states are f\/ixed by the action of an additional discrete symmetry group.)
The choice of codomain plays a role in the analysis in that we will focus on invariant functions
of the form $F = f \circ \iota$, where $f$ is a function on the codomain of $\iota$.
We choose to work with invariant maps with smooth codomains (possibly
with boundary) to avoid possible analytic complications. In particular, we wish to avoid stratum-by-stratum calculations.
The relationship between invariant mappings and quotient spaces has a chicken-and-egg aspect~-- if the quotient space
in a given situation is well-understood and the codomain of a quotient mapping is or can be extended to a manifold, this
is a natural source of an adequate map. However, development of an explicit avatar of an abstract quotient involves
construction of an appropriate invariant mapping; if such an avatar is not known a priori, the construction of a quotient mapping
can involve additional steps beyond those required to construct an adequate map. In particular,  an adequate map is
not required to be surjective.

The map $\iota$ in the last example has a `bonus' invariance, namely invariance under transposition:
$\iota(A^T) = \iota(A)$ for any $A \in {\rm GL}(3)$. $\iota$ is, in fact, invariant with respect to the action
of  the group $\widehat G = \bZ_2 \ltimes ({\rm O}(3) \times {\rm O}(3))$,  with semi-direct product structure
\begin{gather}\label{ad_trans}
(\tau, U, V) (\tau', U', V') =  \begin{cases}
(\tau', U U', V V'), &  \tau = 1, \\
(-\tau', U V', V U') , & \tau = -1.
\end{cases}
\end{gather}
This additional symmetry plays a crucial role in the analysis of the Riemann ellipsoids \cite{Chandrasekhar,FassoLewis, Riemann},
which are discussed in Section~\ref{riemann}.
If we identify ${\rm so}(3)$ with $\R^3$ as usual, and $\fg$ with $\R^3 \times \R^3$, the adjoint action satisf\/ies
\[
\Ad {(\tau, U, V)}(\bseta, \bszeta)  =  \begin{cases}
(U \bseta , V \bszeta), & \tau = 1 ,\\
(U \bszeta , V \bseta), & \tau = -1,
\end{cases}
\]
with associated Lie bracket
\[
[(\bseta, \bszeta), (\tilde \bseta, \tilde \bszeta)] = (\bseta \times \tilde \bseta,\bszeta \times \tilde \bszeta).
\]
A left action of $\widehat G$ on ${\rm SL}(3)$ is determined by
\[
(1, U, V) \cdot A =
U   A V^T
\qquad{\rm and}\qquad
(-1, \idm, \idm) \cdot A = A^T.
\]

Unrecognized or unintentional symmetries of invariant maps can lead to unexpected dege\-ne\-racies of
the dif\/ferentials. If such symmetries are found, it is necessary to decide whether they are, in fact, relevant
to the system or systems at hand; if not, it may be advisable to modify the map $\iota$ so as to break this
extra symmetry.

\section{Relative critical points}
\label{rel_crit_pts}

A simplifying feature of the examples treated in Sections~\ref{lag_top} and \ref{pendulum}
is the commutivity of the symmetry groups. In those examples, the discrete symmetries~-- negation and
ref\/lections across vertical planes~-- determined the nontrivial isotropy subgroups other than the full symmetry group.
In general, locked Lagrangians, energy-momentum functions, etc.  are only invariant with respect to the action
of the isotropy subgroup of the parametrizing algebra element under the adjoint action (or coadjoint action if
parametrized by $\fg^*$).  However, the action of the full
symmetry group typically plays a crucial role in the analysis of the relative critical points. If the (dual) algebra elements
are treated as f\/ixed parameters, invariant function theory cannot be fully utilized.
Hence we extend the group action using the (co)adjoint action
to obtain a fully invariant function on $M \times \fg$ (or $M \times \fg^*$), but only require criticality with respect to the original variables.

Consider a manifold $M$ with a left $G$ action; assume that the action of another group~$\tG$ on~$\fg$ commutes with the adjoint action of
$G$, and that $F: M \times \fg \to \R$ is invariant under the action
\[
(g, h) \cdot (m, \bsxi) = (g \cdot m, h \cdot \Ad g \bsxi)
\]
of $G \times \tG$ on $P = M \times \fg$. ($\fg$ can be replaced with the dual of the algebra $\fg^*$ and the adjoint action replaced with
the coadjoint action throughout.) Consideration of the group $\tG$ is motivated by the importance of time reversal symmetry in many
Lagrangian and Hamiltonian systems; negation of the velocity (or momentum) corresponds to negation of the algebra element when
working with the locked Lagrangian, or with the element of the dual of the algebra when working with the amended potential.

\begin{definition}
$(m, \bsxi)$ is a {\it relative critical point} of $F$ if\/f $m$ is a critical point of the function $F_\bsxi: M \to \R$ given by
$F_\bsxi(m) := F(m, \bsxi)$.
\end{definition}

The following  observation is an immediate consequence of the chain rule, but central to our approach to analyzing relative critical points:
\bigskip

{\it If $F = f \circ \iota$ for some invariant map $\iota: M \times \fg \to \cQ$ and function $f: \cQ \to \R$, then
$(m, \bsxi)$ is a relative critical point of $F$ iff}
\begin{gather}\label{comp_crit_cond}
df(\iota(m, \bsxi)) \in \ker \lsb  d_m^*  \iota_\xi \rsb,    \qquad \mbox{where} \qquad \iota_\xi(m) := \iota(m, \bsxi).
\end{gather}

Clearly, any pair $(m , \bsxi)$ such that $\iota(m , \bsxi)$ is a critical point of $f$ is a relative critical point of~$f \circ \iota$.
To f\/ind any additional relative critical points, we can identify the pairs $(m , \bsxi)$ for which~$d_m^* \iota_\xi$ is not injective, and
then determine the pairs for which~$0 \neq df(\iota(m, \bsxi)) \in \ker \lsb d_m^* \iota_\xi \rsb$.
To illustrate the approach, we begin with a very simple and familiar example, the free rigid body.

\begin{example}[the free rigid body]
The Lagrangian for the free rigid body is obtained by setting $g = 0$ in (\ref{lag_top_lag}).
The Lagrangian is invariant with respect to left multiplication on ${\rm SO}(3)$; the associated inf\/initesimal generator is
$\bsxi_M(A) = \widehat \bsxi A$. If we def\/ine $f: \R^3 \to \R$ by
\[
f(\bsOmega) := \half \bsOmega^T \bbI \bsOmega,
\]
then the locked Lagrangian for the free rigid body satisf\/ies
\[
L_\bsxi(A) = L(\bsxi_M(A)) = L\big(\widehat \bsxi A\big) = L\big(A \widehat {A^T \bsxi} \big) = f\big(A^T \bsxi\big).
\]

Since $\widehat{U \bsxi} = U \widehat \bsxi U^T$ for any $U \in {\rm SO}(3)$ and $\bsxi \in \R^3$,
the adjoint action has the form $\Ad U \bsxi = U \bsxi$ when ${\rm so}(3)$ is identif\/ied with $\R^3$, the map $\iota: {\rm SO}(3) \times \R^3 \to \R^3$
given by $\iota(A, \bsxi) := A^T \bsxi$
is invariant with respect to the action of ${\rm SO}(3)$ on ${\rm SO}(3) \times \R^3$, and  the action of $\bbZ_2$ on
$\R^3$ by negation. The locked Lagrangian satisf\/ies $L_\bsxi = f \circ \iota$.
Mechanically, $\iota$ maps the pair $(A, \bsxi)$ to the body representation of $\bsxi_M(A)$.
If we identify $T_A {\rm SO}(3)$ with $\R^3$ using the left trivialization, then
\[
d_A \iota_\xi(A \widehat \bseta, 0) = - \widehat \bseta A^T \bsxi = \iota(A, \bsxi) \times \bseta = \widehat {\iota(A, \bsxi)} \bseta
\]
implies that $d^*_A \iota_\xi(A) = - \widehat {\iota(A, \bsxi)} $, with kernel
\[
\mbox{ker} [d^*_A \iota_\xi(A)] = \mbox{span} \{\iota(A, \bsxi)\} \qquad \mbox{if} \qquad \bsxi \neq 0.
\]
Combining this with $df(\bsOmega) = \bbI \bsOmega$ yields the relative critical point condition that
\[
\bbI \iota(A, \bsxi) \in \mbox{span} \{\iota(A, \bsxi)\},
\]
which is the familiar relative equilibrium condition that the body angular velocity $\iota(A, \bsxi)$ be an eigenvector of $\bbI$.
\end{example}

The free action and underlying Lie group structure make the rigid body amenable to analysis by a wide variety of approaches (see, e.g.,
\cite{ArnoldMMCM, MarsdenRatiu}).
Note that in the treatment given above, there is no need to use reduction to obtain the Lie--Poisson structure on $\R^3 \approx {\rm so}(3)$. (Of course, the
Lie--Poisson structure for the free rigid body is probably more widely known than the corresponding unreduced
system on $T {\rm SO}(3)$, but the same strategy can be applied in situations where the reduced Poisson structure is not readily available.)

The interaction between $M$ and $\fg$ can be clarif\/ied by choosing manifolds $\cQ_M$, $\cQ_\fg$, and $\cQ_\chi$, a~$G$-adequate map $\phi: M \to \cQ_M$ and $G \times \tG$-adequate map  $\alpha: \fg \to \cQ_\fg $,
and then seeking a map  $\chi: M \times \fg \to \cQ_\chi$ such that the map
$\iota: M \times \fg  \to \cQ := \cQ_M \times \cQ_\chi \times \cQ_\fg$ given by
\begin{gather}\label{iota_triple}
\iota(m, \bsxi) :=  (\phi(m), \chi(m, \bsxi), \alpha(\bsxi))
\end{gather}
is $G \times \tG$-adequate for $M \times \fg$. Note that if points $(m, \bsxi)$ with $G_m \not \subseteq G_\bsxi$ exist, then when constructing an adequate map,
it is essential to choose $\chi$ with suf\/f\/iciently high rank at $(m, \bsxi)$ to compensate for the possible drop in rank of $d_m \phi$
relative to nearby points.
Given a scalar function  $f$ on $\cQ = \cQ_M \times \cQ_\chi \times \cQ_\fg$, def\/ine the partial derivatives
$\frac {\partial f}{\partial \phi}(p, x, a) \in T_p^* \cQ_M$ and
$\frac {\partial f}{\partial \chi}(p, x, a) \in T_x^* \cQ_\chi$ by
\[
\pd f \phi(p, x, a)(\delta p) = d f(p, x, a)(\delta p, 0, 0)
\qquad {\rm and}\qquad
\pd f \chi(p, x, a)(\delta x) = d f(p, x, a)(0, \delta x,  0).
\]

\begin{proposition}\label{rel_crit_conds}
$(m, \bsxi)$ is a relative critical point of $F = f \circ \iota$ for $\iota$ satisfying \eqref{iota_triple} iff
\begin{gather}\label{comp_crit}
d_m^* \phi \lp \frac {\partial f}{\partial \phi}(\iota(m, \bsxi))   \rp + d_m^* \chi_\bsxi \lp \frac {\partial f}{\partial \chi}(\iota(m, \bsxi))   \rp = 0,
\end{gather}
where $\chi_\bsxi(m) := \chi(m, \bsxi)$ for all $m \in M$. In particular, $(m, \bsxi)$ is a relative critical point if
\begin{gather}\label{comp_crit_inter}
\pd f \phi(\iota(m, \bsxi)) \in \ker \lsb d_m^* \phi \rsb
\qquad \text{and}\qquad
\pd f \chi(\iota(m, \bsxi))  \in \ker \lsb d_m^*\chi_\bsxi \rsb.
\end{gather}
If $\rang{d_m^*  \phi} \cap \rang{d_m^*  \chi_\bsxi}$ is trivial, then \eqref{comp_crit_inter} is a necessary condition for relative
criticality.
\end{proposition}

\begin{proof}
The f\/irst assertion follows immediately from the identity
\begin{gather*}
d(f \circ \iota)(m, \xi)(\delta m, 0)  =  \pd f \phi(\iota(m, \bsxi))(d_m \phi(\delta m)) + \pd f \chi(\iota(m, \bsxi))(d_m \chi_\xi(\delta m)) \\
\hphantom{d(f \circ \iota)(m, \xi)(\delta m, 0)}{}
 =  \lp d_m^* \phi \lp \frac {\partial f}{\partial \phi}(\iota(m, \bsxi))    \rp + d_m^* \chi_\bsxi \lp \frac {\partial f}{\partial \chi}(\iota(m, \bsxi))   \rp \rp (\delta m)
\end{gather*}
for all $(m, \bsxi)$ and $\delta m \in T_m M$.

If (\ref{comp_crit}) is satisf\/ied, but the two terms in the sum are nonzero, then
\[
d_m^* \phi \lp \frac {\partial f}{\partial \phi}(\iota(m, \bsxi))   \rp = -  d_m^* \chi_\bsxi \lp \frac {\partial f}{\partial \chi}(\iota(m, \bsxi))   \rp \neq 0;
\]
hence $\rang{d_m^*  \phi} \cap \rang{d_m^*  \chi_\bsxi}$ must be nontrivial.
\end{proof}

Given an invariant mapping $\iota: P \to \cQ$, (\ref{comp_crit_cond}) can be used to identify functions $f: \cQ \to \R$ such that a given pair $(m, \bsxi)$ is
a relative critical point of $f \circ \iota$. This can be used in model development, e.g., if it is known a priori that a given state $m$ should maintain
steady motion with generator $\bsxi$, then Lagrangian or Hamiltonian models of the system can be designed for which $(m, \bsxi)$ is a
relative critical point of the associated invariant function.
This approach can be incorporated into data f\/itting schemes as follows: Given a sequence of
experimental data, motion capture measurements, or numerical simulations, and a choice of conf\/iguration manifolds,  a best
f\/it conf\/iguration can be assigned to each data point using (nonlinear) least squares methods. However, if it is known a priori
that the evolution of the system closely approximates steady motion, it my be advantageous to f\/ind not the closest conf\/igurations
in the model conf\/iguration space, but the closest conf\/igurations consistent with steady motion. If some known or postulated features of
the system impose nontrivial conditions on the derivatives of any function that might be used to determine the dynamics in a
Hamiltonian or Lagrangian model, then data collected from an approximately steady motion can be projected onto states compatible
with the criticality condition~(\ref{comp_crit_cond}).

\subsection{Rigid conditions and rigid-internal decompositions}

A striking feature of the free rigid body is the coincidence of the conf\/iguration manifold and symmetry group;
the codomain of the invariant map $\iota$ is the Lie algebra of the symmetry group. In general, not all relevant
information about a  system is captured by the symmetry group. However, consideration of
the criticality conditions arising from variations tangent to the group orbits in $M \times \fg$ can be an important step
in the derivation and interpretation of the full criticality conditions. The coupling between the given action on $M$
and the coadjoint action on $\fg$ (or $\fg^*$) induced by the invariance of the function or map often yields signif\/icant
computational simplif\/ications.

Given a manifold $N \neq \R$ and mapping $\Psi: M \times \fg \to N$, let $\pd \Psi \bsxi(m, \bsxi): \fg \to T_{\Psi(m, \xi)} N$
denote the linear map
\[
\pd \Psi \bsxi(m, \bsxi)(\bseta) := d_{(m, \bsxi)} \Psi(0, \bseta).
\]
Given a scalar function $F: M \times \fg \to N$, def\/ine $\pd F\bsxi(m, \bsxi) \in \fg^*$ by $\pd F \bsxi(m, \bsxi)\cdot \bseta := dF(m, \bsxi)(0, \bseta)$
for all $\bseta \in \fg$.

An invariant function $F:P \to \R$ satisf\/ies
\[
dF(m, \bsxi)(\bseta_M(m), \ad \bseta \bsxi) = 0
\]
for any $\bseta \in \fg$; hence $dF_\bsxi(m)|_{\fg \cdot m} = 0$ if\/f
\begin{gather}\label{rigid}
\ads \bsxi \pd F \bsxi (m, \bsxi) = 0.
\end{gather}
We will refer to (\ref{rigid}) as the {\em rigid condition}.

\begin{example}[energy-momentum method]
 Consider a symplectic manifold $M$ with $G$-invariant Hamiltonian $H$ and
$\Ads{}$-equivariant momentum map $\bJ: M \to \fg^*$. Set $F(m, \bsxi) = H(m) - j_\bsxi(m)$, where
$j_\bsxi: M \to \R$ is given by $j_\bsxi(m) = J(m) \cdot \bsxi$. The function $F$ is linear in $\bsxi$, with
$\pd F \bsxi(m, \bsxi) = - J(m)$. Hence the rigid condition is $\ads \bsxi J(m) = 0$
(see, e.g., \cite{LS89b, LS89a,LMSP89,SPM90}).
\end{example}

\begin{example}[locked Lagrangian]Given a $G$-invariant Lagrangian $L$ on $TQ$, with a lifted $G$ action,
def\/ine the {\it locked momentum map} by $\bbI_ \bsxi := J \circ \bsxi_Q$, where $J$ is
the momentum map associated to the lifted action; equivalently,
\[
\bbI_\bsxi(q) \cdot \bszeta = \dep {L_{\bsxi + \eps   \bszeta}(q)}.
\]
The rigid condition for $F(q, \bsxi) = L_\bsxi(q) = L(\bsxi_Q(q))$ is $\ads \bsxi \bbI_ \bsxi(q) = 0$  (see \cite{Lewis92}).
\end{example}

If $F = f \circ \iota$, then (\ref{rigid}) is satisf\/ied if\/f
\begin{gather*}%\label{iota_rigid}
d f(\iota(m, \bsxi)) \in \ker \lsb \lp  \pd \iota \bsxi (m, \bsxi) \circ \ad \bsxi \rp^* \rsb.
\end{gather*}

If a pair $(m, \bsxi)$ satisf\/ies the rigid condition for a map $F$, the relative criticality of the pair can be determined by
considering the restriction of $dF_\bsxi(m)$ to a  complement $\vint$ to $\fg \cdot m$ in~$T_m M$.
Criticality conditions complementary to the rigid condition are called {\it internal conditions}. A~rigid-internal
decomposition for the energy-momentum method on cotangent bundles of principal bundles was introduced in \cite{LMSP89},
and subsequently adapted to the reduced energy-momentum method in \cite{SLM91}. An analogous
decomposition for locked Lagrangians, including systems with nonfree actions, was developed in \cite{Lewis92}; this decomposition
can be used to construct a basis with respect to which the linearized dynamics of Lagrangian systems at relative equilibria
block diagonalize. A convenient choice of complement can facilitate the analysis, but convenience is typically context-dependent.
In some approaches, the tangent space of  an appropriate (local) submanifold  can be taken as $\vint$
in other situations, $\vint$ is specif\/ied via a system of linear equations on $T_m M$. For example, the decompositions in
\cite{Lewis92, LMSP89, SLM91} are determined by  the linearization of the rigid condition with respect to~$m$ and block diagonalize
the second variation of the relevant function. Thus these choices of internal variations are well-suited for stability and bifurcation analyses.

The following corollary to Proposition~\ref{rel_crit_conds} establishes the relative criticality conditions with respect to a
rigid-internal decomposition.

\begin{corollary}
\label{rel_crit_conds_cor}
If $F = f \circ \iota$ for some invariant map $\iota$ and scalar function $f$,  and $\vint \subseteq T_m M$ satisfies
$T_m M = \vint + \fg \cdot m$,
then $(m, \bsxi)$ is a relative critical point iff
\begin{gather}\label{rig_int}
d f(\iota(m, \bsxi)) \in \cK_{\rm rig}(m, \bsxi; \iota) \cap \ker \lsb \lp d\iota_\xi (m)|_{\vint} \rp^* \rsb,
 \end{gather}
 where
 \begin{gather}\label{K_rig}
 \cK_{\rm rig}(m, \bsxi; \iota):=  \ker \lsb \lp \pd \iota \bsxi (m, \bsxi) \circ \ad \bsxi \rp^* \rsb.
\end{gather}
 If, in addition, $\iota$ has the form \eqref{iota_triple}, then \eqref{rig_int} is equivalent to
 \begin{gather}
  \pd f \chi (\iota(m, \bsxi)) \in \cK_{\rm rig}(m, \bsxi; \chi), \nonumber \\
 \lp d_m \phi|_{\vint} \rp^*\lp \pd f \phi (\iota(m, \bsxi))\rp + \lp d_m \chi_\xi|_{\vint}\rp^* \lp \pd f \chi (\iota(m, \bsxi))\rp = 0 ,
\label{chi_rigid}
\end{gather}
for $\cK_{\rm rig}(m, \bsxi; \chi)$ defined analogously to \eqref{K_rig}.

If $\mbox{\rm range} \big[ \lp d_m \phi|_{\vint} \rp^* \big]$ and $\mbox{\rm range} \big[ \lp d_m \chi_\xi|_{\vint}\rp^* |_{\cK_{\rm rig}(m, \bsxi; \chi)}\big]$
have trivial intersection, then \eqref{chi_rigid} is satisfied iff
\begin{gather*}
\pd f \phi (\iota(m, \bsxi)) \in \ker \lsb \lp d_m \phi|_{\vint} \rp^* \rsb
\qquad\!\! \text{and}\!\!\qquad
\pd f \chi (\iota(m, \bsxi)) \in \cK_{\rm rig}(m, \bsxi; \chi) \cap \ker \lsb \lp d_m \chi_\xi|_{\vint} \rp^* \rsb.
\end{gather*}
 \end{corollary}

\begin{proof}
The proof follows directly from Proposition~\ref{rel_crit_conds} and the chain rule.
\end{proof}

Note that changing the dependence of $f$ on $\phi$, or changing the map $\phi$ itself, does not
change the rigid condition. Hence this condition can be particularly useful when studying the relative
critical points of families of functions with diverse purely $M$-dependent components. For example,
changing the potential energy in an invariant natural mechanical system does not alter the associated rigid condition.

If $M$ is a principal $G$-bundle with  projection $\pi: M \to B$, then
given a neighborhood $\cU$ of~$\pi(m)$ and a local section $\bssigma: \cU \to \cV \subseteq M$, we can take
$\vint = d_{\pi(m)} \bssigma(T_{\pi(m)} \cU)$; $(m, \bsxi)$ is a~relative critical point of $F$ if\/f $(m, \bsxi)$
satisf\/ies the rigid conditions and $\pi(m)$ is a critical point of the function
$b \mapsto F(\bssigma(b), \bsxi)$. If the action is not free, the image of a local section can often be replaced by a~slice
(see, e.g., \cite{CLOR02} and references therein).
A submanifold $S$ of $M$ is  a {\it  slice} at $m \in S$ if \vspace{-.05in}
\begin{itemize}\itemsep=0pt
\item
$S$ is invariant under $G_m$;
\item
 the union $G \cdot S$ of all orbits intersecting $S$ is an open neighborhood of the orbit $G \cdot m$;
\item
there is an equivariant mapping $\pi: G \cdot S \to G \cdot m$ such that $\pi|_{G \cdot m}$ is the identity on
$G \cdot m$ and $\pi^{-1}(m) = S$.
\end{itemize}

If $F = f \circ \iota$ for some invariant map $\iota$ and scalar function $f$, and $S$ is a slice through a pair $(m, \bsxi)$
satisfying the rigid conditions, then we can take $\vint = T_m S$. Hence
internal criticality conditions can be obtained by replacing the restriction of $d_m \iota_\bsxi$ to $\vint$ with the
linearization $d_m (\iota_\xi|_S)$ of the restriction of $\iota_\bsxi$ to $S$ in~(\ref{comp_crit_cond}).

\section{The symmetries of the Riemann ellipsoids}\label{riemann}

The Riemann ellipsoids provide a historically signif\/icant example with nonabelian symmetry, and hence nontrivial rigid component, and
diverse assortment of isotropy subgroups. The Riemann ellipsoids are steady motions of an incompressible inviscid f\/luid moving under
a self-gravitational potential that are generated by one parameter subgroups of the action of ${\rm SO}(3)$ on the left and the group of volume
preserving dif\/feomorphisms on the right~\cite{Riemann}. Rigidly rotating ellipsoids were used as qualitative models of spinning f\/luid masses
by Newton, Maclaurin, and
Jacobi \cite{Chandrasekhar}. Dirichlet \cite{Dirichlet} rigorously proved that af\/f\/ine motions of an ellipsoidal f\/luid mass form an
invariant subsystem of the Euler equations
for an incompressible inviscid f\/luid under the inf\/luence of a~self-gravitational potential. (See \cite{Chandrasekhar} for a~survey of the history of the
Riemann ellipsoids and a~modern derivation using moments. See, e.g.,  \cite{FassoLewis,Lewis93, ROSD} for analyses using
geometric methods\footnote{Riemann utilized energy criteria when considering the stability and bifurcation of the relative equilibria;
precursors of both the reduced energy-momentum and energy-Casimir methods can be found in \cite{Riemann}.
\cite{Chandrasekhar} carried out a traditional linear stability analysis of the relative equilibria and found that some families of equilibria
that failed to meet the energy criteria considered by Riemann were nonetheless linearly stable.

Recent studies of the Riemann ellipsoids using modern geometric
methods include the following: \cite{Lewis93} includes an analysis of the bifurcation of families of triaxial relative equilibria from the Maclaurin ellipsoids,
which are axisymmetric and in steady rotation about the axis of symmetry; this paper uses the locked Lagrangian.  In \cite{FassoLewis}
Nekhoroshev stability methods are used to investigate the long-time behavior of those ellipsoids that are linearly stable but fail to meet the `natural' energy method
nonlinear stability criterion. This paper includes an overview of the linear and nonlinear stability of the ellipsoids; these
results are derived and expressed using modern geometric constructs. (The linear and nonlinear stability classif\/ications in \cite{FassoLewis} are exact,
but the Nekhoroshev stability analysis involves numerical approximations in the computation of the Birkhof\/f normal forms.) More recently,
\cite{ROSD} uses the augmented potential and techniques analogous to those used in \cite{Lewis93} to analyze
the nonlinear stability and linear (in)stability results for axisymmetric relative equilibria in steady rotation about a principal axis of the f\/luid mass.
}.)

The  af\/f\/ine motions of an ellipsoidal f\/luid mass can be characterized as a
Lagrangian~-- or canonical Hamiltonian~-- system on ${\rm SL}(3)$; the right action of volume preserving dif\/feomorphisms on the original
conf\/iguration manifold determines an action of ${\rm SO}(3)$ by right multiplication.
The af\/f\/ine f\/luid model has conf\/iguration manifold ${\rm SL}(3)$ and symmetry group $\bZ_2 \ltimes ({\rm SO}(3) \times {\rm SO}(3))$ with
multiplication (\ref{ad_trans}).
The analysis of the relative equilibria of this system traces the development of symmetry methods in conservative systems.
Maclaurin  analyzed the rigid rotations of axisymmetric f\/luid masses about the axis of symmetry;
Jacobi  extended these results to triaxial f\/luid masses~\cite{Chandrasekhar}.
Dedekind \cite{Dedekind} observed that rigidly rotating steady motions, generated by
left multiplication, are transformed by transposition to internal swirling steady motions, in which the region occupied by the f\/luid
remains unchanged. This transposition invariance played a crucial role
in the evolution of symmetry-based analyses of dynamical systems, extending techniques developed for
readily visualizable, intuitively clear spatial rotations to novel group actions~-- the full body of results for the left action immediately yielded
a corresponding body of results for the action of ${\rm SO}(3)$ by right multiplication associated to particle relabeling invariance. Riemann \cite{Riemann}
built on Dedekind's insight, classifying all relative equilibria~-- including large families of quasi-periodic motions~-- generated by combinations of left and right multiplication.

We now generalize some of the classical results for the Riemann ellipsoids by determining relative criticality conditions for
the action of $G = \bZ_2 \ltimes ({\rm SO}(3) \times {\rm SO}(3))$ on ${\rm SL}(3)$ obtained by restricting the action of
$\bZ_2 \ltimes ({\rm O}(3) \times {\rm O}(3))$ on ${\rm GL}(3)$ described in the last example in Proposition~\ref{adequate_examples} and
the remarks following the proposition. We construct an invariant map $\iota$ on ${\rm SL}(3) \times \R^3 \times \R^3$ associated to that action,
identifying the Lie algebra of ${\rm SO}(3)$ with $\R^3$ as usual, and develop relative critical point conditions for functions of the form $f \circ \iota$.
Given the relative complexity of this system, we do not rigorously investigate the adequacy of this map here, noting only that it satisf\/ies
(\ref{adequate_cond}) at generic points within some key isotropy classes. We consider two classes of points $(A, \bseta, \bszeta) \in P
= {\rm SL}(3) \times \R^3 \times \R^3$ studied by Riemann. S type points are those for which both components of the algebra element are
parallel to a common axis determined by $A$; specif\/ically, an S type point is one for which there is a singular value decomposition
$A = U \mbox{diag}[\ba] V^T$, $U, V \in {\rm O}(3)$, such that $\bseta = \eta \, U \be_j$ and $\bszeta = \zeta \, V \be_j$ for some $j \in \{1, 2, 3\}$ and
$\eta, \zeta \in \R$.
Riemann's type I, II, and III ellipsoids all have coplanar algebra elements: there is a singular value decomposition $A = U \mbox{diag}[\ba] V^T$
such that $U^T \bseta$ and $V^T \bszeta$ belong to a common coordinate plane.
As we shall see, these states have nontrivial  isotropy that reduces the number of criticality conditions to be satisf\/ied. This analysis generalizes
one of the central components of Riemann's treatment of the af\/f\/ine liquid drop model, extending some of his results for locked Lagrangians and
amended potentials determined by the ${\rm SO}(3) \times {\rm SO}(3)$ action and the self-gravitational potential to a much more general class of invariant
functions. This generalization sheds light on the crucial role played by symmetry in Riemann's original analysis, as well as more recent
treatments.

\begin{proposition}
\label{phi_alpha}
The map $\phi: {\rm SL}(3) \to \lp \R^+ \rp^2$ given by
\begin{gather}\label{sl3_iota}
\phi(A) :=  \big(\norm{A}^2, \norm{A^{-1}}^2 \big)
\end{gather}
is adequate for the action of $G = \bZ_2 \ltimes ({\rm SO}(3) \times {\rm SO}(3))$ on ${\rm SL}(3)$, and is equivariant with respect to the $\bZ_2$ actions
\[
A \mapsto A^{-1}
\qquad \text{and}\qquad
(x, y) \mapsto (y, x).
\]
on ${\rm SL}(3)$ and $(\R^+)^2$.

The map
\begin{gather}\label{alpha_def}
\alpha(\bseta, \bszeta) := \lp \norm{\bseta}^2 + \norm{\bszeta}^2, \big( \norm{\bseta}^2 - \norm{\bszeta}^2 \big)^2 \rp
\end{gather}
is adequate for the adjoint action on $\fg \approx \R^3 \times \R^3$.
\end{proposition}

\begin{proof}
If $A \in {\rm SL}(3)$, then
\[
\norm{A}^2 = \mbox{trace} \, \big(A^T A\big)
\qquad {\rm and}\qquad
 \norm{A^{-1}}^2 =  \half \lp (\mbox{trace} \, (A^T A))^2 - \mbox{trace} \lp (A^T A)^2 \rp \rp.
 \]
Hence $\phi(A)$ can be regarded as the restriction of the invariant
map def\/ined on ${\rm GL}(3)$ in Proposition~\ref{adequate_examples}, with the third component~-- which is identically equal to one on ${\rm SL}(3)$~-- discarded.
(\ref{sl3_iota})~satisf\/ies
\[
\mbox{rank} [d_A \phi] + 1  =  \mbox{\# of distinct singular values of $A$}.
\]
The verif\/ication of adequacy of $\phi$ follows the argument given for ${\rm GL}(3)$ in Proposition~\ref{adequate_examples}.

The linearization $d_{(\eta, \zeta)} \alpha$ satisf\/ies
\[
\half   d_{(\eta, \zeta)} \alpha =   \begin{pmatrix}  \bseta^T & \bszeta^T \\
c   \bseta^T & - c   \bszeta^T \end{pmatrix},
\qquad \mbox{where} \qquad c = 2 \big( \norm{\bseta}^2 - \norm{\bszeta}^2 \big).
\]
If $0 \neq \norm{\bseta} \neq \norm{\bszeta} \neq 0$, then $d_{(\eta, \zeta)} \alpha$ has rank two and
$\fg \cdot (\bseta, \bszeta) = \{ \bu \times \bseta, \bv \times \bszeta): \bu, \bv \in \R^3 \}$ is four-dimensional.

If $0 \neq \norm{\bseta} = \norm{\bszeta}$, then $\bszeta = V \bseta$ for some $V \in {\rm SO}(3)$ and
$(\bseta, \bszeta)$ is f\/ixed by $(-1, V^T, V)$. The linearization of the action of $(-1, V^T,V)$ satisf\/ies
\[
d_{(\eta, \zeta)} \Phi_{(-1, V^T,V)}(\bx, \by) = \big(V^T \by , V \bx \big),
\]
and hence
\[
\mbox{range} \big[ d_{(\eta, \zeta)} \Phi_{(-1, V^T,V)} - \idm \big] = \big\{ (\bx, - V \bx) : \bx \in\R^3 \big\} =
\ker [ d_{(\eta, V \eta)} \alpha ].
\]
Since $\fg \cdot (\bseta, V \bseta) = \{ (\bx \times \bseta, \by \times V \bseta) :\bx, \by \in \R^3 \}$,
\[
\fg \cdot (\bseta, V \bseta) + \mbox{range} \big[ d_{(\eta, \zeta)} \Phi_{(-1, V^T,V)} - \idm \big]
\]
is f\/ive-dimensional and (\ref{adequate_cond}) is satisf\/ied at $(\bseta, V \bseta)$.

If $0 = \norm{\bseta} \neq \norm{\bszeta}$, then
$\mbox{range} \lsb d_{(\eta, \zeta)} \alpha \rsb = \mbox{span} \{ (1, - c) \}$ is one-dimensional and
$\fg \cdot (0, \bszeta) = \{ 0, \bv \times \bszeta):  \bv \in \R^3 \}$. The isotropy subgroup
of $(0, \bszeta)$ is $G_{(0, \bszeta)} = \{1 \} \times {\rm SO}(3) \times \{\mbox{rotations about}\ \bszeta \}$, and
\begin{gather*}
\cup_{U \in {\rm SO}(3), \theta \in S^1} \mbox{range} \lsb d_{(0, \zeta)} \Phi_{(1, U, \exp(\theta \bszeta))} - \idm  \rsb
\supseteq \cup_{U \in {\rm SO}(3)}\big\{(U   \bx - \bx, 0) : \bx \in \R^3 \big\}
= \R^3 \!\times \{0 \}.
\end{gather*}
Hence $\fg \cdot (0, \bszeta) +
\cup_{U \in {\rm SO}(3), \theta \in S^1} \mbox{range} \lsb d_{(0, \zeta)} \Phi_{(1, U, \exp(\theta \bszeta))} - \idm  \rsb$
is at least f\/ive-dimensional, so (\ref{adequate_cond}) is satisf\/ied at $(0, \bszeta)$. The verif\/ication of
(\ref{adequate_cond}) when $0 = \norm{\bszeta} \neq \norm{\bseta}$ is entirely analogous.

Finally, $G_{(0, 0)} = G$ and $\cup_{U, V \in {\rm SO}(3)} \mbox{range} \lsb d_{(0, 0)} \Phi_{(1, U, V)} - \idm  \rsb = \R^3 \times \R^3$
implies that $\alpha$ is adequate.
\end{proof}

Note that the map
$(\bseta, \bszeta) \mapsto  \big( \norm{\bseta}^2, \norm{\bszeta}^2 \big)$ is not invariant under the adjoint action
(\ref{ad_trans}) of group elements involving transposition, since that action exchanges $\bseta$ and $\bszeta$.

We now choose a map $\chi$ depending nontrivially on both $M$ and $\fg$.
At generic points $d_{(A, \bseta, \bszeta)} (\phi \times \alpha)$ has rank four; however, if $A$ has exactly two distinct
singular values~-- and hence $d_A \phi$ is rank one~-- but the isotropy subgroup of $(\bseta, \bszeta)$ has trivial intersection
with that of~$A$, then $\phi \times \chi \times \alpha$ cannot be adequate at $(A, \bseta, \bszeta)$ unless $\chi$ has
codomain of dimension greater than four.
In light of the equivariance of $\phi$ with respect to inversion of $A$, and the fact that the entries of $A^{-1}$ are polynomials
in the entries of $A$ for $A \in {\rm SL}(3)$, we will use inversion in the construction of the invariant map $\iota$.
This equivariance simplif\/ies the analysis of some important classes of relative critical points.
Def\/ine scalar functions $K$ and $C$, and vector-valued function~$\chi$ on $P = {\rm SL}(3) \times \R^3 \times \R^3$ as follows:
\begin{gather*}
K(A, \bseta, \bszeta) := \half \norm{(\bseta, \bszeta)_M(A)}^2 =  \half \big\|\hat \bseta A - A \hat \bszeta\big\|^2,
\\
C(A, \bseta, \bszeta) := \lp \bseta \times A A^T \bseta \rp^T A \bszeta + \lp \bszeta \times A^T A \bszeta \rp^T A^T \bseta,
\end{gather*}
and
\begin{gather}\label{chi_def}
\chi(A, \bseta, \bszeta) :=
\lp  K(A, \bseta, \bszeta), K(A^{-T}, \bseta, \bszeta), \bseta^T A \bszeta,  \bseta^T A^{-T} \bszeta,  C(A, \bseta, \bszeta) \rp.
\end{gather}

\begin{proposition}
The map $\iota = \phi \times \chi \times \alpha: P \to \R^2 \times \R^5 \times \R^2$ determined by \eqref{sl3_iota}, \eqref{chi_def}, and \eqref{alpha_def}
is invariant under the action of $G = \bZ_2 \ltimes ({\rm SO}(3) \times {\rm SO}(3))$ on $P$ and
equivariant with respect to the action of $\bZ_2^2$ on $P$  by
\begin{gather}\label{extra_sym}
(-1, 1) \cdot (A, \bseta, \bszeta) = (A, -\bseta, - \bszeta) \qquad \text{and}\qquad
(1, -1) \cdot (A, \bseta, \bszeta) = (A^{-T}, \bseta, \bszeta),
\end{gather}
and the linear $\bZ_2^2$ action on $\R^2 \times \R^5 \times \R^2 \approx \R^9$ such that $(-1, 1)$ negates the seventh component, leaving the other components unchanged,
and $(1, -1)$ pairwise exchanges the first six components and negates the seventh component.
\end{proposition}

\begin{proof}
 Invariance of $\phi$ and $\alpha$ with respect to the action of $G$ was demonstrated in Proposition~\ref{phi_alpha}. Invariance of
the matrix norm with respect to transposition implies that $\phi(A^{-T}) = \phi(A^{-1})$, which has components equal to the
permuted components of $\phi(A)$. $\alpha$ is clearly invariant with respect to the $\bZ_2^2$ action.

The identity $\widehat{U \bsxi} = U \widehat \bsxi U^T$ for any $U \in {\rm SO}(3)$ and $\bsxi \in \R^3$ determines the invariance of
$\chi$ with respect to elements of $G$ of the form $(1, U, V)$.
\[
\hat \bszeta A^T - A^T \hat \bseta = \big(\hat \bseta A - A \hat \bszeta\big)^T
\]
shows that $\chi$ is also invariant with respect to the action $(-1, \idm, \idm)
 \cdot (A, \bseta, \bszeta) = (A^T,  \bszeta, \bseta)$. Hence~$\chi$ is invariant with respect to the action of $G$.

The function $K$ and weighted inner products are invariant with respect to negation of the algebra element, while $C$ is equivariant with respect to this action and the combination of inversion and transposition of $A$, specif\/ically
\[
C(A, -\bseta, - \bszeta) = - C(A, \bseta, \bszeta) = C\big(A^{-T}, \bseta, \bszeta\big).
\]

That the action (\ref{extra_sym}) of $(1, -1)$ pairwise exchanges the f\/irst  six  entries of $\iota$ and leaves the last two unchanged is obvious.
The anti-symmetry of $C$ with respect to the action of $(1, -1)$
follows from the invariance of the triple product with respect to multiplication of all three terms by an element of ${\rm SL}(3)$. Specif\/ically,
if $\vert \bx \ \ \by \ \ \bz \vert$ denotes the determinant of the matrix with columns $\bx$, $\by$, $\bz$, then
\begin{gather*}
\lp \bseta \times A A^T \bseta \rp^T A \bszeta  =  \left \vert \bseta \ \ A A^T \bseta \ \ A \bszeta \right \vert
 =  \det  \big( A A^T \big) \big \vert  \lp A A^T \rp^{-1}  \bseta \ \ \bseta \ \ A^{-T}  \bszeta \big \vert \\
\hphantom{\lp \bseta \times A A^T \bseta \rp^T A \bszeta }{}
=  - \big \vert \bseta \ \ A^{-T}   A^{-1} \bseta \ \ A^{-T}   \bszeta \big  \vert
 =  - \big( \bseta \times   A^{-T}   \lp A^{-T} \rp^T \bseta \big)^T   A^{-T}  \bszeta
\end{gather*}
and, analogously,
\begin{gather*}
\lp \bszeta \times A^T A \bszeta \rp^T A^T \bseta
 =
\det  \lp A^T A \rp \big \vert  \lp A^T  A \rp^{-1}  \bszeta \ \ \bszeta \ \ A^{-1} \bseta \big \vert \\
\hphantom{\lp \bszeta \times A^T A \bszeta \rp^T A^T \bseta}{}
 =
- \big( \bszeta \times \lp A^{-T}\rp^T   A^{-T}   \bseta \big)^T  \lp A^{-T} \rp^T  \bszeta
\end{gather*}
imply $C(A, \bseta, \bszeta) + C(A^{-T}, \bseta, \bszeta) = 0$.
\end{proof}

\subsection{Coplanar states}

The time reversal symmetry of the Lagrangian for the af\/f\/ine f\/luid system leads to invariance of the locked Lagrangian under negation of the algebra element. The interaction between this action and the $G$ action on $P$ plays an essential role in Riemann's analysis of the af\/f\/ine liquid drop system \cite{Riemann}. The invariance of the map $\iota$ with respect to the $G$  action and equivariance with respect to the $\bZ_2^2$ action~(\ref{extra_sym})
determine much of the structure of the dif\/ferential of $\iota$ at points with nontrivial isotropy with respect to these actions.

We will say that $(A, \bseta, \bszeta)$ is {\it coplanar} if $\bseta$ and $\bszeta$ both lie in plane determined by
$A$, i.e.\  if there exist $U, V \in {\rm SO}(3)$ such that $U A V^T$ is diagonal, and $\mbox{span} \{ U \bseta, V \bszeta \}$ is the Euclidean orthogonal
complement of one of the Euclidean axes. Riemann's  types I, II, and III classes of relative equilibria are coplanar \cite{Riemann}.
The coplanar states are f\/ixed by the two element subgroup of the enlarged group
$G \times \bZ_2^2$ generated by ref\/lection across the plane containing both components of the algebra element.  Specif\/ically,
given $\ba \in \R^3$ with $a_1 a_2 a_3 = 1$, $\bseta, \bszeta \in \R^3$, and $U, V \in {\rm SO}(3)$, if we def\/ine
\begin{gather*}
\varphi(\ba, \bseta, \bszeta, U, V) := \big(U \dia V^T, U   \bseta, V   \bszeta\big)
\qquad\!\!\! {\rm and}\qquad\!\!
g_j(U, V) := \big(1, U R_j U^T, V  R_j V^T\big)\! \in \! G,
\end{gather*}
where $R_j$ denotes rotation about $\be_j$ through $\pi$, then
\[
g_j(U, V)  \cdot \varphi(\ba, \bseta, \bszeta, U, V) = \phi(\ba, \bseta, \bszeta, U  R_j, V  R_j).
\]
In particular, if  $\bseta$ and $\bszeta$ are both orthogonal to $\be_j$, then
$(g_j(U, V), (-1, 1))$ f\/ixes $\varphi(\ba, \bseta, \bszeta, U, V)$.

The equivariance of $C$ with respect to negation of the algebra elements implies that
coplanar points are zeroes of $C$: if $p = \varphi(A, \bseta, \bszeta, U, V)$ is coplanar and
$\psi : P \to P$ denotes the action of $(g_j(U, V), (-1, 1))$, then $C \circ \psi = - C$ implies
$- C(p) = C(\psi(p)) = C(p)$.

Many of the following calculations were carried out using {\sl Mathematica}'s symbolic capabilities; hence the results are stated without detailed
proofs. We have identif\/ied subspaces that lie in the relevant kernels for all points within the considered classes, but
satisfaction of appropriate algebraic relations between components of $A$, $\bseta$, and $\bszeta$ can led to a drop in rank of $d_A \iota_{(\eta, \zeta)}$
relative to generic points within the same class. We do not carry out a detailed analysis of all the possibilities here; rather, we focus on behavior that
is clearly determined by the symmetries outlined above.

For notational convenience, assume that $p = \varphi(\dia, \bseta, \bszeta)$, with $\bseta_3 = \bszeta_3 = 0$; the general coplanar case follows from invariance.
Let $(r_j, \theta_j/2)$ be polar coordinates of $(\bseta_j, \bszeta_j)$, $j = 1, 2$.
(The motivation for this angle convention will become clear.) We f\/irst consider the rigid conditions:
If $a_1 \neq a_2$, then the space
$\cK_{\rm rig}(p)$ of solutions of the rigid condition (\ref{chi_rigid}) satisf\/ies
\begin{gather}\label{kr_def}
\cK_{\rm rig}(p) \supseteq \mbox{span}\lcb \bv(\ba, \theta_1, \theta_2),  \widetilde \bv(\ba, \theta_1, \theta_2) \rcb,
\end{gather}
where
\begin{gather}\label{kr_bv}
\bv(\ba, \theta_1, \theta_2) :=(\cos \theta_2 - \cos \theta_1)\lp \half \be_1 + 2   \be_4 \rp - \bw(1, 2) + \bw(2,1),
\\
\bw(i, j) := \cos \theta_i \sin \theta_j \lp a_j   \be_3 - \frac {a_i}{a_3}   \be_4 \rp,\nonumber
\end{gather}
and $\widetilde \bv(\ba, \theta_1, \theta_2)$ is constructed by replacing $a_i$ with $\frac 1 {a_i}$, $i = 1, 2, 3$, in $\bv(\ba, \theta_1, \theta_2)$
and exchanging the f\/irst and second, and the third and fourth, components of the resulting expression. If $\dia$ has three distinct singular values, then $\fg \cdot \dia$ is six-dimensional and we can take the submanifold of diagonal elements of
${\rm SL}(3)$ as a slice through $\dia$.
There exist matrix-valued functions $N, \tilde N: \R^3  \to \R^{2 \times 5}$ such that the internal
criticality condition is satisf\/ied if\/f
\begin{gather*}%\label{riem_inter_tri}
\pd f \phi(\iota(p)) = \big( r_1^2 \big( N(\ba)  + \tilde N(\ba) \sin \theta_1\big)
 + 1 \leftrightarrow 2 \big)  \pd f \chi(\iota(p)),
\end{gather*}
where by $1 \leftrightarrow 2$ we mean that the f\/irst and second components of $\ba$
are exchanged, $r_1$ is replaced by $r_2$, and $\theta_1$ replaced by $\theta_2$. Thus in this case, given
$\pd f \chi(\iota(p)) \in \cK_{\rm rig}(p)$, there is a unique value of $\pd f \phi(\iota(p))$ compatible with relative
criticality of $p$ for $f \circ \iota$.

We now consider axisymmetric coplanar ellipsoids: If $a_1 = a_2 = a$, and hence $a_3 = a^{-2}$, then
\begin{gather*}%\label{kr_def_sym}
\cK_{\rm rig}(p) \supseteq \mbox{span}\lcb \be_1 + 4   \be_4, \be_2 + 4   \be_3, \be_3 - a^2 \be_4 \rcb.
\end{gather*}
If we take
\[
\tbe_i := (\be_i,  \bzero, \bzero), \quad i = 1, 2,
\qquad {\rm and}\qquad
\abe_i := (\bzero, \be_i, \bzero), \quad i = 1, \ldots, 5,
\]
as basis vectors for $T_{\iota(p)}^* \cQ$, then the internal criticality condition is satisf\/ied if $df(\iota(p))$ lies in
\[
df(\iota(p)) \in \mbox{ker}[d_A^*\iota(p)]
= \mbox{span} \lcb a^2  \tbe_1 - \tbe_2, w(a, \bseta, \bszeta)  \tbe_1 + \abe_1 + 4   \abe_4 + a^4 (\abe_2 + 4   \abe_3)  \rcb,
\]
where
\[
w(a, \bseta, \bszeta) := - \lp 1+ a^6 \rp \big(  \norm{\bseta}^2 + \norm{\bszeta}^2 \big).
\]
Thus
$\mbox{range} \big[ d_{\dia}^*  \phi|_S \big] \cap \mbox{range} \big[ \lp d_p^*  \chi_S \rp|_{\cK_{\rm rig}(p)} \big]$ is one-dimensional for generic values of $a$;
hence there are relative critical points with nonzero partial derivative with respect to $\chi$.
Note that in both cases~-- two and three distinct singular values~-- $\mbox{ker}[d_A^*\iota(p)]$ is at least two-dimensional, although the spaces of partial derivatives satisfying the rigid criticality conditions are of dif\/ferent dimensions. (The nullity of $d_A^*\iota(p)$ can be greater if appropriate algebraic
conditions relating the components of $A$, $\bseta$, and $\bszeta$ are satisf\/ied.)

\begin{remark}
To recapture the rigid conditions for the locked Lagrangian for the actual Riemann ellipsoid problem, with algebra-dependent component
$K(A, \bseta, \bszeta)$, we f\/ind conditions under which $\be_1 \in \cK_{\rm rig}(p)$. Considering (\ref{kr_def}) and (\ref{kr_bv}), we
see that solutions of the equations
\begin{gather*}
a_1 \cos \theta_2 \sin \theta_1  =  a_2 \cos \theta_1 \sin \theta_2, \\
a_2 \cos \theta_2 \sin \theta_1 - a_1 \cos \theta_1 \sin \theta_2  =  2   a_3 (\cos \theta_2 - \cos \theta_1) \neq 0
\end{gather*}
determine states satisfying the rigid condition for the locked Lagrangian. If appropriate non-degeneracy conditions on
$\theta_1$ and $\theta_2$ are satisf\/ied, these equations have the unique solution
\[
\frac {a_1}{a_3} = \frac {2   \cos \theta_1 \sin \theta_2}{\cos \theta_1 + \cos \theta_2},
\qquad
\frac {a_2}{a_3} = \frac {2   \cos \theta_2 \sin \theta_1}{\cos \theta_1 + \cos \theta_2},
\qquad {\rm and}\qquad
a_3^3 = \frac {(\cos \theta_1 + \cos \theta_2)^2}{\sin 2 \theta_1   \sin 2 \theta_2}.
\]
\end{remark}

\subsection{S type states}

We brief\/ly consider the S type ellipsoids, for which $\bseta$ and $\bszeta$ both lie on the same Euclidean coordinate axis
after some change of basis that diagonalizes $A$.
If $p = \varphi(\ba, \eta   \be_j, \zeta   \be_j, U, V)$ for $\ba \in \cA$, $\eta, \zeta \in \R$, $j \in \{1, 2, 3\}$, and $U, V \in {\rm SO}(3)$,
then $p$ satisf\/ies the rigid condition, and hence is a relative critical point if it satisf\/ies the internal conditions.

If the absolute values of the entries of $\ba$ are distinct, then
\[
\mbox{ker}[d^*_p \iota] \supseteq
\mbox{span}\lcb \bv_{12}\lp \smallfrac 1{a_3} \rp, \bv_{21}(a_3), \sin\theta  (a_3 \tbe_1 - \tbe_2) + 4   c_1 \abe_3,  \abe_3 + a_3^3 \abe_4, \abe_5 \rcb,
\]
where $(\eta, \zeta) = \lp r   \cos \smallfrac \theta 2, r   \sin \smallfrac \theta 2 \rp$ and $\bv_{ij}(u) :=  c_1 \abe_j + c_2(u) \tbe_i + c_3(u) \tbe_j$
for
\[
 c_1 := \frac {\lp a_1^2 - a_3^2 \rp \lp a_3^2 - a_2^2 \rp}{r^2 a_3},
\qquad c_2(u) := u^3 - \frac {a_1}{a_2} - \frac {a_2}{a_1} + \sin \theta,
\qquad
c_3(u) := \frac 1 u - u^2 \sin \theta.
\]

If $p$ has continuous isotropy component $S^1$, say $p = \varphi \lp \lp a, a, a^{-2} \rp, \eta   \be_3, \zeta   \be_3, U, V \rp$, $|a| \neq 1$, then
$d_p^* \iota$ is rank one, with
\[
\mbox{ker}[d_p^* \iota] \supseteq
\mbox{span}\lcb (\eta - \zeta)^2 \tbe_1 +c   \abe_1, \eta   \zeta   \tbe_1 + 2   c    \abe_4, a^2 \tbe_1 - \tbe_2, \abe_1 + a^4 \abe_2, a^4 \abe_3 + \abe_4 , \abe_5 \rcb,
\]
where $c := a^{-6} - 1 $.

\section{Conclusion}

The relative critical point method is closely related to several well-established approaches to the analysis of relative equilibria,
as well as more general critical point problems.
One crucial distinction between this formulation and others, e.g.\ the energy-Casimir and (reduced) energy-momentum
methods, is the absence of any geometric structures associated to Hamiltonian or Lagrangian mechanics. This not only
allows the application of the method in more general settings, but disengages the purely symmetry-driven aspects of the
analysis of relative equilibria of conservative systems from the Poisson, symplectic, or variational structures linking the
functions in question to the dynamics. Thus one gains insight into the roles of the key features of the
system at each stage of the formulation and solution of the problem.

Replacing families of functions parametrized by the algebra (or its dual) of the symmetry group with an invariant function
facilitates formulation of the criticality conditions as linear conditions on the partial derivatives of a function on the codomain
of an invariant map, as in Proposition~\ref{rel_crit_conds} and Corollary~\ref{rel_crit_conds_cor}. The resulting systems of
equations are well-suited both for the classif\/ication of relative critical points and for the identif\/ication
of invariant functions for which a given state is a relative critical point. No quotient
manifolds or orbifolds need be determined; changes in isotropy class typically lead to changes in dimension of the subspaces in which
admissible partial derivatives must lie, but do not result in singularities. Discrete symmetries are easily managed using invariant
maps, and play an important role in the analysis of relative critical points. The relevant maps are invariant with respect to the full
symmetry group~-- there is no need to restrict attention to the isotropy subgroup of a f\/ixed algebra element.
Thus the relative critical point approach can be particularly
advantageous in the analysis of systems possessing a rich symmetry group.

\subsection*{Acknowledgments}

The author is indebted to the referees for their many valuable suggestions and corrections. Their
insightful contributions greatly improved this work.

\pdfbookmark[1]{References}{ref}
\LastPageEnding

\end{document}